\documentclass[12pt,reqno]{amsart}

\usepackage{amsthm,amsfonts,amssymb,euscript}
\numberwithin{equation}{section} 
\newcommand{\eps}{\epsilon}
\newcommand{\slim}{\operatorname{s-lim}}
\begin{document}

\title[Agmon--Kato--Kuroda theorems]{Agmon--Kato--Kuroda theorems for a large class of perturbations}

\author{Alexandru D. Ionescu}
\address{Department of Mathematics, University of Wisconsin--Madison, Madison, WI 53706, U.S.A.} \email{ionescu@math.wisc.edu}
\author{Wilhelm Schlag}
\address{Department of Mathematics, The University of Chicago, 5734 South University Avenue, Chicago, IL 60637, U.S.A.}
\email{schlag@math.uchicago.edu}
\thanks{Both authors were supported in part by NSF grants and Alfred P. Sloan research fellowships.}

\begin{abstract} We prove asymptotic completeness for operators of the form $H=-\Delta+L$
on $L^2(\mathbb{R}^d)$, $d\ge2$, where $L$ is an {\it{admissible
perturbation}}. Our class of admissible perturbations contains
multiplication operators defined by real-valued potentials $V\in
L^q(\mathbb{R}^d)$, $q\in[d/2,(d+1)/2]$ (if $d=2$ then we require
$q\in(1,3/2]$), as well as real-valued potentials $V$ satisfying a
global Kato condition. The class of admissible perturbations also
contains first order differential operators of the form
$\vec{a}\cdot\nabla- \nabla\cdot\overline{\vec{a}}$ for suitable vector
potentials $a$. Our main technical statement is a new limiting
absorption principle which we prove using techniques from harmonic
analysis related to the Stein-Tomas restriction theorem.
\end{abstract}

\maketitle

\section{Introduction}

One of the basic problems of quantum mechanics is to determine
the spectrum and the spectral types of the self-adjoint operator
\[ H = -\Delta +L\]
on $L^2(\mathbb{R}^d)$, where $L$ is a suitable perturbation. A
minimal requirement for self-adjointness is that $L$ is symmetric.
Given the self-adjoint operator $H$, let
$\sigma_{\mathrm{ac}}(H)$, $\sigma_{\mathrm{sc}}(H)$, and
$\sigma_{\mathrm{pp}}(H)$, denote its absolutely continuous
spectrum, singular continuous spectrum, and pure point spectrum
respectively. Let $H= \int \lambda \, E(d\lambda)$ denote the
spectral resolution of $H$. It is well-known that there is a
Lebesgue decomposition
\[ E= E_{\rm ac}+E_{\rm sc}+E_{\rm pp},\]
where the terms on the right-hand side are projection valued
measures. The ranges of $E_{\rm ac}(\mathbb{R})$, $E_{\rm
sc}(\mathbb{R})$, and $E_{\rm pp}(\mathbb{R})$ are orthogonal and
are typically denoted by $L^2_{\rm ac}, L^2_{\rm sc}$, and
$L^2_{\rm pp}$, respectively. These subspaces are referred to as
the absolutely continuous, singular continuous, and pure point
subspaces, respectively. Physically, it is most relevant to
determine which of these is nonzero. This is related to the
long-time behavior of the evolution $e^{-itH}$. Indeed, any $f\in
L^2(\mathbb{R}^d)$ with $E_{\rm pp}f=f$ does not propagate,
whereas $(E_{\rm ac}+E_{\rm sc})f=f$ leads to transport (see the
RAGE theorem in~\cite{CFKS}).

A much-studied case of perturbations are those defined by multiplication
with suitable potentials $V$. For example, if for some $\varepsilon>0$
\begin{equation}
\label{eq:Vdec}
 \sup_{x\in\mathbb{R}^d} (1+|x|)^{1+\varepsilon}|V(x)| < \infty,
\end{equation}
then a classical theorem of S.~Agmon \cite{Ag} (which applies to all dimensions $d\ge1$),
combined with T.~Kato's theorem \cite{Ka} on
absence of eigenvalues in $(0,\infty)$ for such $V$,
states there is {\em asymptotic completeness} in this case. In dynamical terms, this refers
to the fact that for any $f\in L^2(\mathbb{R}^d)$ there is $f_0\in L^2(\mathbb{R}^d)$
so that
\[ e^{-itH}f = \sum_{j} e^{-it\lambda_j} P_{j} f + e^{-itH_0} f_0 + o_{L^2}(1) \]
as $t\to\infty$. Here $H_0=-\Delta$, $\lambda_j\le0$ are the
eigenvalues of $H$, and $P_{j}$ are the orthogonal projections
onto the associated eigenspaces. In spectral terms, this means
that $E_{\rm sc}=0$ and that the wave operators
\[ \Omega^{\pm}:= \slim_{t\to\mp\infty} e^{iHt} e^{-itH_0} \]
exist and are complete, i.e., they are surjective onto the absolutely continuous
spectral subspace $L^2_{\rm ac}$ of~$H$.
S.~Agmon's work was the culmination of a series of partial results for which we refer to~\cite{Ag}
and M.~Reed, B.~Simon~\cite{RS4}. In particular, Agmon deduced the existence and completeness of the
wave operators from the {\em limiting absorption principle}
\begin{equation}
\label{eq:lim1}
\sup_{\lambda\ge\lambda_0} \sup_{\epsilon>0}
\|(H-(\lambda+i\epsilon))^{-1}\|_{L^{2,\sigma}\to L^{2,-\sigma}} \le C(V,\lambda_0),
\end{equation}
$\lambda_0>0$ and $\sigma>1/2$,
via Kato's smoothing theory, see~\cite{kato}. Here
\[ L^{2,\sigma}(\mathbb{R}^d):= \big\{ f\::\: (1+|x|)^\sigma f(x)\in L^{2}(\mathbb{R}^d) \big\}.\]
We remark that~\eqref{eq:lim1} immediately leads to the fact that
\[ E_{\rm sc}(\mathbb{R}^+)=0\]
 because of the density of $L^{2,\sigma}$ in~$L^2$,
see Theorem~XIII.20 in~\cite{RS4}. A recent example of A.~Kiselev
\cite{kis} (in $d=1$) shows that S.~Agmon's theorem is essentially
sharp as far as the decay of $V$ is concerned. For a recent review
of much of what is known about the spectral theory of decaying
potentials we refer to S.~Denisov, A.~Kiselev's
survey~\cite{DenKis}.

The optimality of~\eqref{eq:Vdec} is related to the optimality of
$\sigma>1/2$ in the limiting absorption principle~\eqref{eq:lim1}.
When $V\equiv 0$, the limiting absorption principle
\eqref{eq:lim1} is intimately connected to basic {\em restriction
theorems} for the Fourier transform. The relevant restriction
theorem in this case is the bound
\[ \|\hat{f}\|_{L^2(\mathbb{S}^{d-1})} \le C\|f\|_{L^{2,\sigma}(\mathbb{R}^d)} \]
with $\sigma>1/2$, known as the {\em trace-lemma}.

The trace lemma applies to the restriction of the
Fourier transform to any compact hypersurface. In particular, it does not use the
fact that the Gaussian curvature of the sphere does not vanish.
In contrast, the well-known Stein-Tomas restriction theorem asserts that
\[ \|\hat{f}\|_{L^2(\mathbb{S}^{d-1})} \le C\|f\|_{L^{p_d}(\mathbb{R}^d)} \]
where $p_d=(2d+2)/(d+3)$ and $d\ge2$ (see \cite{To}). This is an
optimal bound in the sense that it fails for any $p>p_d$.
Moreover, it fails for surfaces with one vanishing principal
curvature. It is natural to ask what kind of Agmon-type theorem or
limiting absorption principle results from using the Stein-Tomas
theorem rather than the much simpler trace lemma. This issue was
addressed by M.~Goldberg and the second author~\cite{golsch} who obtained
the bound
\begin{equation}
\label{eq:V34}
\sup_{0<\epsilon<1,\;\lambda\ge\lambda_0}
\Big\|(-\Delta+V - (\lambda^2+i\epsilon))^{-1}\Big\|_{L^{4/3}\to L^4} \le C(\lambda_0,V)\;\lambda^{-1/2}
\end{equation}
provided $V\in L^{3/2}(\mathbb{R}^3)\cap
L^{3/2+\delta}(\mathbb{R}^3)$, $\delta>0$. In particular, the
spectrum of $-\Delta+V$ is purely absolutely continuous on
$(0,\infty)$ for such~$V$. This result depended on the recent unique
continuation theorem of the first author and D.~Jerison~\cite{IoJe},
who established the absence of imbedded point spectrum for $H$ under
the condition $V\in L^{3/2}(\mathbb{R}^3)$ (with suitable analogues
in all dimensions $d\ge2$). Because of its dependence on a strong
unique continuation result at infinity, the approach
of~\cite{golsch} was rather limited. In particular, it applied only
to potentials $V \in L^{3/2}(\mathbb{R}^3)\cap
L^{3/2+\delta}(\mathbb{R}^3)$, $\delta>0$. Moreover,
in~\cite{golsch} no unconditional statement could be made about
absence of singular continuous spectrum for $V\in
L^p(\mathbb{R}^3)$, $3/2\le p\le 2$.

The goal of this paper is to prove an Agmon-type theorem for a
much larger class of perturbations without relying on any unique
continuation theorem at infinity. Similar to~\cite{Ag}, we will
prove a suitable limiting absorption bound. However, extensive use
is made of bounds on oscillatory integrals in the spirit of the
Stein-Tomas restriction theorem and related bounds for
Bochner-Riesz means, see \cite[Chapter IX]{St}. We now describe
our results in more detail. Our main theorem is Theorem~\ref{main1}.

We assume from now on that the dimension $d$ is $\geq 2$. We
define the sets $D_j=\{x\in\mathbb{R}^d:|x|\in[2^{j-1},2^{j}]\}$,
$j\geq 1$, and $D_0=\{x\in\mathbb{R}^d:|x|\in[0,1]\}$. Following
the notation in \cite[Chapter XIV]{Ho}, we also define the
following Banach spaces of functions on $\mathbb{R}^d$, $d\geq 2$:
\begin{equation*}
\begin{split}
&B=\Big\{f:\mathbb{R}^d\to\mathbb{C}:||f||_B:=\sum_{j=0}^\infty
2^{j/2}||f||_{L^2(D_j)}<\infty\Big\},\\
&B^\ast=\Big\{u:\mathbb{R}^d\to\mathbb{C}:||u||_{B^\ast}:=\sup_{j\geq
0
}2^{-j/2}||u||_{L^2(D_j)}<\infty\Big\}.
\end{split}
\end{equation*}
The spaces $B$ and $B^\ast$ are related to the sharp form of the trace lemma
\begin{equation}\label{shtr}
\mathcal{F}:B\to L^2(\mathbb{S}^{d-1})\text{ and }\mathcal{F}^{-1}:L^2(\mathbb{S}^{d-1})\to B^\ast,
\end{equation}
as bounded operators (see \cite[Theorem 7.1.26]{Ho}).

Let $\mathcal{S}(\mathbb{R}^d)$ denote the space of Schwartz
functions on $\mathbb{R}^d$ and $\mathcal{S}'(\mathbb{R}^d)$ the
space of distributions. For any $\alpha\in\mathbb{C}$ let
$S_\alpha:\mathcal{S}'(\mathbb{R}^d)\to\mathcal{S}'(\mathbb{R}^d)$
denote the operator defined by the Fourier multiplier
$\xi\to(1+|\xi|^2)^{\alpha/2}$, i.e.,
$S_\alpha=(1-\Delta)^{\alpha/2}$. For $1<p<\infty$ and
$\alpha\in\mathbb{R}$ we define the standard Sobolev spaces
\begin{equation*}
W^{\alpha,p}=\{u\in\mathcal{S}'(\mathbb{R}^d):S_\alpha u\in L^p\}\text{ with }||u||_{W^{\alpha,p}}:=||S_\alpha u||_{L^p}.
\end{equation*}
Let $p_d=(2d+2)/(d+3)$ and $p'_d=(2d+2)/(d-1)$ denote the
Stein-Tomas restriction exponents. Let $S_1(B)$ denote the image
of $B$ under $S_1$, and $S_{-1}(B^*)$ the image of $B^*$ under
$S_{-1}$. The main Banach spaces we use in this paper are
\begin{equation*}
X:=W^{-1/(d+1),p_d}+S_1(B)\text { with }||f||_{X}:=\inf_{f_1+f_2=f}||S_{-1/(d+1)}f_1||_{L^{p_d}}+||S_{-1}f_2||_{B},
\end{equation*}
and
\begin{equation*}
X^\ast:=W^{1/(d+1),p'_d}\cap S_{-1}(B^\ast)\text { with }||u||_{X^\ast}:=\max(||S_{1/(d+1)}u||_{L^{p'_d}},||S_{1}u||_{B^\ast}).
\end{equation*}
Clearly, $X$ is a space of distributions and $X^\ast\subseteq
W^{1,2}_{\mathrm{loc}}$. To motivate these definitions, we notice
first that
\begin{equation}\label{tb1}
\mathcal{F}:X\to L^2(\mathbb{S}^{d-1})\text{ and }\mathcal{F}^{-1}:L^2(\mathbb{S}^{d-1})\to X^\ast,
\end{equation}
as bounded operators, which is a consequence of the Stein-Tomas
restriction theorem and \eqref{shtr}. Moreover
\begin{equation}\label{tb2}
X\hookrightarrow W^{-1,2}\text{ and }W^{1,2}\hookrightarrow X^\ast,
\end{equation}
which follows from the Sobolev imbedding theorem (this explains
the choice of the exponent $1/(d+1)$ in the definition of $X$ and
$X^\ast$) and the imbedding $B\hookrightarrow L^2\hookrightarrow
B^\ast$. Finally, for more general theorems, we would like to have
the space $X$ as large as possible and the space $X^\ast$ as small
as possible, subject to \eqref{tb1} and \eqref{tb2}. Our first
theorem is a uniform bound for the free resolvent $R_0(z)$,
$z\in\mathbb{C}\setminus[0,\infty)$. The operator $R_0(z)$ is
defined on $\mathcal{S}'(\mathbb{R}^d)$ by the Fourier multiplier
$\xi\to(|\xi|^2-z)^{-1}$.

\newtheorem{main2}{Theorem}[section]
\begin{main2}\label{main2}
Assume that $\delta\in(0,1]$. Then
\begin{equation}\label{bu1}
\sup_{|\lambda|\in[\delta,\delta^{-1}],\epsilon\in[-1,1]\setminus\{0\}}||R_0(\lambda+i\epsilon)||_{X\to X^\ast}\leq
C_\delta<\infty,
\end{equation}
where $C_\delta$ is a (finite) constant that depends only on
$\delta$ and the dimension $d$.
\end{main2}

The main point of Theorem~\ref{main2} is the uniformity of the bound
\eqref{bu1} as $\epsilon\to 0$. In contrast, the bound in the
stronger (elliptic) imbedding $R_0(\lambda+i\epsilon):W^{-1,2}\to
W^{1,2}$ blows up as $\epsilon\to 0$ if $\lambda>0$. The proof of
Theorem \ref{main2} is essentially known through the work of
L.~H\"ormander, C.~Kenig, A.~Ruiz, C.~Sogge, and L.~Vega,
see~\cite{Ho}, \cite{KeRuSo}, \cite{RuVe}; we collect the necessary
bounds in Section~\ref{free}.

We also prove a weighted estimate. For $N\geq 0$,
$\gamma\in(0,1]$, and $x\in\mathbb{R}^d$, we define the weight
\begin{equation}\label{tb3}
\mu_{N,\gamma}(x)=\frac{(1+|x|^2)^N}{(1+\gamma |x|^2)^N}.
\end{equation}

\newtheorem{main3}[main2]{Theorem}
\begin{main3}\label{main3}
Assume that $\delta\in(0,1]$. Then
\begin{equation}\label{bo1}
||\mu_{N,\gamma}u||_{X^\ast}\leq
C_{N,\delta}||\mu_{N,\gamma}(\Delta+\lambda)u||_{X}
\end{equation}
for any $\lambda\in\mathbb{R}$ with
$|\lambda|\in[\delta,\delta^{-1}]$, and any $u\in X^\ast$ with the
property that
\begin{equation}\label{bo100}
\lim_{R\to\infty}R^{-1}\int_{R\leq |x|\leq 2R}|u|^2\,dx=0.
\end{equation}
The constant $C_{N,\delta}$ depends only on $N$,
$\delta$, and the dimension $d$.
\end{main3}

We remark that the condition \eqref{bo100} is necessary: let
\begin{equation*}
u(x)=\int_{\mathbb{S}^{d-1}}e^{-ix\cdot\xi}\,d\xi.
\end{equation*}
Then $u\in X^\ast$, however $(\Delta+1)u\equiv 0$.
Theorem~\ref{main3} plays a key role in the bootstrap argument in
the proof of our main Theorem~\ref{main1} below (see
Lemma~\ref{lemma9.1}). We emphasize that the constant in
\eqref{bo1} is allowed to depend on the parameter $N$, but not on
$\gamma\in(0,1]$.

For functions $u,f\in\mathcal{S}(\mathbb{R}^d)$, we define $\langle
u,f\rangle:=\int_{\mathbb{R}^d}u\overline{f}\,dx$. Clearly,
$\langle u,f\rangle=\langle S_\alpha(u),S_{-\alpha}(f)\rangle$ for
any $\alpha\in\mathbb{R}$. By a slight abuse of notation, we
extend the definition of $\langle.,.\rangle$ to pairs
\begin{equation*}
(u,f)\in \mathcal{S}'(\mathbb{R}^d)\times\mathcal{S}(\mathbb{R}^d)\cup L^2\times L^2\cup W^{-1,2}\times W^{1,2}\cup X\times X^\ast.
\end{equation*}
We have
\begin{equation}\label{tb5}
|\langle u,f\rangle|\leq\min(||u||_{L^2}||f||_{L^2},||u||_{W^{-1,2}}||f||_{W^{1,2}},||u||_{X}||f||_{X^\ast}).
\end{equation}
Also, it follows easily from the definitions of the spaces $X$ and $X^\ast$ that
\begin{equation}\label{tb511}
||f||_{X}\leq C\sup_{\phi\in\mathcal{S}(\mathbb{R}^d),\,||\phi||_{X^\ast}=1}|\langle f,\phi\rangle|\text{ and }||u||_{X^\ast}\leq C\sup_{\phi\in\mathcal{S}(\mathbb{R}^d),\,||\phi||_{X}=1}|\langle u,\phi\rangle|.
\end{equation}

{\bf{Definition:}} {\em  Let $\mathcal{L}(X^\ast,X)$ denote the space of
bounded operators from $X^\ast$ to $X$. We say that $L$ is an
{\em{admissible perturbation}} if:

(1) $L\in\mathcal{L}(X^\ast,X)$ and
\begin{equation}\label{tb61}
\langle L\phi,\psi\rangle=\overline{\langle L\psi,\phi\rangle}
\end{equation}
for any $\phi,\psi\in \mathcal{S}(\mathbb{R}^d)$ (i.e.\ $L$ is
symmetric).

(2) For any $\varepsilon>0$ and $N\geq 0$ there are
$A_{N,\varepsilon},R_{N,\varepsilon}\in[1,\infty)$ such that
\begin{equation}\label{tb6}
||\mu_{N,\gamma}Lu||_X\leq
\varepsilon||\mu_{N,\gamma}u||_{X^\ast}+A_{N,\varepsilon}||u\mathbf{1}_{\{|x|\leq
R_{N,\varepsilon}\}}||_{L^2}
\end{equation}
for any $u\in X^\ast$ and any $\gamma\in(0,1]$, where
$\mathbf{1}_E$ denotes the characteristic function of the set $E$.

(3) There is an integer $J\geq 1$ and operators $A_j, B_j$ in
$\mathcal{L}(X^\ast,L^2)$ for $1\le j\le J$ such that
\begin{equation}
\label{eq:Lsplit} \langle L \phi,\psi \rangle = \sum_{j=1}^J
\langle B_j \phi, A_j \psi \rangle\text{ for any } \phi, \psi \in
X^*.
\end{equation}
Moreover, considered as (unbounded) operators on $L^2$, $A_j,B_j$
are closed on some domains satisfying\footnote{These inclusions
are natural since $W^{1,2}\subset X^*$}
\[ \mathrm{Domain}(A_j)\supseteq W^{1,2}({\mathbb R}^d),\quad
\mathrm{Domain}(B_j)\supseteq W^{1,2}({\mathbb R}^d)\] for all
$1\le j\le J$. }
\medskip

Our main goal is to prove an Agmon-type theorem for admissible
perturbations. Before formulating our main theorem, we remark that
condition~(1) in the definition of admissible perturbations is
essential for our arguments. Condition~(2) is somewhat technical
and is related to our use of Theorem~\ref{main3}. Variations (and
improvements) of condition~(2) are possible. Condition~(3) is the
usual condition which arises in Kato's smoothing
theory~\cite{kato}, and is needed in order to study the wave
operators which intertwine $-\Delta$ and $-\Delta+L$.

In view of \eqref{tb2},
\begin{equation*}
H:=-\Delta+L:W^{1,2}\to W^{-1,2}
\end{equation*}
as a bounded operator, for any admissible perturbation $L$. Our
main theorem is the following:

\newtheorem{main1}[main2]{Theorem}
\begin{main1}\label{main1}
Assume that $L$ is an admissible perturbation. Then the following properties hold:

(a) The operator $H=-\Delta+L$ defines a closed, self-adjoint
operator on
\begin{equation}\label{gu1}
\mathrm{Domain}(H)=\{u\in W^{1,2}(\mathbb{R}^d):Hu\in L^2(\mathbb{R}^d)\}.
\end{equation}
In addition, $\mathrm{Domain}(H)$ is dense in $L^2(\mathbb{R}^d)$
and $H$ is bounded from below on $\mathrm{Domain}(H)$.

(b) The set of nonzero eigenvalues ${\mathcal
E}=\sigma_{\mathrm{pp}}\setminus\{0\}$ of $H$ is discrete in
$\mathbb{R}\setminus\{0\}$, i.e., ${\mathcal E}\cap I$ is finite
for any compact set $I\subset\mathbb{R}\setminus\{0\}$. Moreover,
each eigenvalue in $\mathcal{E}$ has finite multiplicity.

(c) Any eigenfunction $u$ of $H$ with eigenvalue $\lambda\ne 0$ is
rapidly decreasing, i.e., for any integer $N\geq 0$,
\begin{equation}\label{bi22}
(1+|x|^2)^Nu\in W^{1,2}(\mathbb{R}^d).
\end{equation}

(d) Let $I\subset (\mathbb{R}\setminus\{0\})\setminus{\mathcal E}$
be compact. Then
\begin{equation}\label{bi3}
\sup_{\lambda\in I,\epsilon\in[-1,1]\setminus
0}||R_L(\lambda+i\epsilon)||_{X\to X^\ast}\leq
C(L,I)<\infty,
\end{equation}
where $R_L(\lambda+i\epsilon)$ denotes the resolvent of $H$ at
$\lambda+i\varepsilon$, and $C(L,I)$ is a constant that
depends on the interval $I$, the perturbation $L$, and the dimension
$d$. Thus the spectrum of the
operator $H$ is purely absolutely continuous on $I$.

(e) $\sigma_{\mathrm{sc}}(H)=\emptyset$ and
$\sigma_{\mathrm{ac}}(H)= [0,\infty)$.

(f) The wave operators $\Omega^{\pm}(H,H_0)$ exist and are
complete, where $H_0=-\Delta$.
\end{main1}

We notice the similarity of Theorem~\ref{main1} with the
Agmon-Kato-Kuroda theorem, see Theorem~XIII.33 in~\cite{RS4}. The
main novelty in our theorem is that it applies to a much larger class of perturbations.
To provide examples of admissible perturbations we define the Banach space
\begin{equation*}
Y=\Big\{V:\mathbb{R}^d\to\mathbb{C}:||V||_Y:=\sum_{j=0}^\infty
2^{j}||f||_{L^\infty(D_j)}<\infty\Big\},
\end{equation*}
where the sets $D_j$ are as in the definitions of the spaces $B$ and $B^\ast$.
For $\delta\in(0,1/2]$ we define the kernels
\begin{equation*}
K_{d,\delta}(x)=\mathbf{1}_{\{|x|\leq\delta\}}
\begin{cases}
|x|^{-(d-2)}\quad&\text{ if }d\geq 3;\\
\log(1/|x|)\quad&\text{ if }d=2.
\end{cases}
\end{equation*}
For any exponent $q\in[1,\infty)$ and measurable function $f$ let
\begin{equation*}
M_q(f)(x)=\Big[\int_{|y|\leq 1/2}|f(x+y)|^{q}\,dy\Big]^{1/q}.
\end{equation*}
Clearly, $M_{q}(f)(x)\leq CM_{q'}(f)(x)$ if $1\leq q\leq
q'\leq\infty$. Also, $||M_q(f)||_{L^{p'}(D_j)}\leq
C||M_q(f)||_{L^p(\widetilde{D}_j)}$ if $1\leq p\leq p'\leq\infty$,
where $\widetilde{D}_j=D_{j-1}\cup D_j\cup D_{j+1}$ if $j\geq 1$
and $\widetilde{D}_0=D_{0}\cup D_1$ (the last inequality is easy
to prove for $p'=p$ and $p'=\infty$, thus for $p'\in[p,\infty]$ by
interpolation). We fix $q_0=d/2$ if $d\geq 3$ and $q_0>1$ if
$d=2$.

\newtheorem{prop16}[main2]{Proposition}
\begin{prop16}\label{prop16}
The following are examples of admissible perturbations:

(a) Multiplication operators defined by real-valued potentials
$V$ with the property that
\begin{equation}\label{tb99}
M_{q_0}(V)\in L^{(d+1)/2},
\end{equation}
or
\begin{equation}\label{tb991}
M_{q_0}(V)\in Y,
\end{equation}
or
\begin{equation}\label{tb100}
\lim_{\delta\to 0}||\,|V|\ast K_{d,\delta}||_Y=0.
\end{equation}

(b) First order differential operators of the form
\begin{equation*}
\vec{a}\cdot\nabla-\nabla\cdot\overline{\vec{a}},
\end{equation*}
defined by vector-valued potentials $\vec{a}:\mathbb{R}^d\to\mathbb{C}^d$ with the property that
\begin{equation}\label{tb98}
\Big[\sum_{j=0}^\infty\Big(2^{j/2}||M_{2q_0}(|\vec{a}|)||_{L^{d+1}(D_j)}\Big)^{p(d)}\Big]^{1/p(d)}<\infty,
\end{equation}
or
\begin{equation}\label{tb981}
M_{2q_0}(|\vec{a}|)\in Y,
\end{equation}
or
\begin{equation}\label{tb101}
\lim_{\delta\to 0}||\,[|\vec{a}|^2\ast K_{d,\delta}]^{1/2}||_Y=0.
\end{equation}

(c) Any finite linear combination of admissible perturbations with real coefficients.
\end{prop16}

We remark that the exponent $(d+1)/2$ in \eqref{tb99} is optimal for Theorem~\ref{main1} to hold.
This is due to a recent example by the first author
 and D.~Jerison~\cite{IoJe} of a potential $V\in L^p$, for
all $p>(d+1)/2$, such that $H=-\Delta+V$ has slowly decaying
eigenfunctions (and positive eigenvalues). We emphasize that this
example is not related to the local singularities of $V$. In fact,
$V$ is a smooth, real-valued function with oscillations and
asymptotic behavior
\[ |V(x)|\approx (1+|x_1|+|x'|^2)^{-1}.\]
The main issue here is that the potential behaves differently
along different directions. It remains to be seen if such examples
can lead to dense point spectrum or even imbedded singular
continuous spectrum as well. It is possible that the transition
point for singular continuous spectrum occurs at larger values
of~$q$, for example at $q=d$ (Coulomb case)\footnote{This
possibility was communicated  to the second author by Barry Simon,
who believes that it should be $q=d$.}. The same remark applies to
first order perturbations defined by vector potentials $a$ as in
\eqref{tb98}. The restriction on the exponent $q_0$ is needed to
define the operator $H$ as a self-adjoint operator on its domain.

In some cases of admissible perturbations we can add the natural conclusion
\begin{equation}\label{tb200}
\sigma_{pp}\subseteq(-\infty,0].
\end{equation}
For potentials $V\in L^{d/2}(\mathbb{R}^d)$, $d\geq 3$ (a
restriction stronger than \eqref{tb99}), this follows from
\cite[Theorem 2.1]{IoJe}. H.~Koch and D.~Tataru~\cite{kochtat} have
recently proved the absence of positive eigenvalues for potentials
$V$ that satisfy conditions similar to~\eqref{tb99}.

We also allow perturbations given by multiplication with
potentials in the global Kato class described in \eqref{tb100}.
For comparison, the {\it{local Kato class}} (cf. \cite{si}) is
defined by the condition
\begin{equation*}
\lim_{\delta\to 0}||\,|V|\ast K_{d,\delta}||_{L^\infty}=0.
\end{equation*}
We remark that the condition \eqref{tb100} is more general than S.
Agmon's condition (1.3) in \cite{Ag}
\[ \sup_{x\in\mathbb{R}^d} \Big[(1+|x|)^{2+2\epsilon} \int_{|y-x|\le 1}
|V(y)|^2 |y-x|^{-d+\mu} \, dy\Big]<\infty.
\]
for some $\epsilon>0$ and $0<\mu<4$. This is easy to see, using
the Cauchy--Schwartz inequality and the fact that
$(1+|x|)^{-1-\epsilon}\in Y$.
For operators defined by potentials
$V$ as in \eqref{tb100}, the conclusion~\eqref{tb200} is not
known; even the easier question of absence of compactly supported
eigenfunctions for such potentials is not settled (see
\cite[p.~519]{si}). Our proof of Theorem~\ref{main1} is
independent of the validity of~\eqref{tb200}.

We have the following corollary of Theorem~\ref{main1}, which
relies on the well-known connection between smoothing bounds for
the resolvent as in~\eqref{bi3} and time-dependent smoothing
bounds. This connection is given by T.~Kato's theory~\cite{kato}.

\newtheorem{mai2}[main2]{Corollary}
\begin{mai2}\label{mai2}
Let $L$ be an admissible perturbation, and let $\mathcal{E}$ be as
in Theorem~\ref{main1}. Then for any compact set
$I\subseteq(\mathbb{R}\setminus\{0\})\setminus\mathcal{E}$, there
exists a constant $C(I,L)$ such that
\begin{equation}\label{mai9}
\begin{split}
||S_{1/(d+1)}[e^{itH}E(I)f]||_{L^{p'_d}_xL^2_t}+||S_{1}[e^{itH}E(I)f]||_{B^\ast_xL^2_t}\le C(I,L) \|f\|_2,
\end{split}
\end{equation}
for any $f\in L^2$, where $E(I)$
denotes the spectral projection onto the interval $I$ associated
with $H$.
\end{mai2}

Because of the spectral projection $E(I)$ with a compact $I$ the
smoothing effect in Corollary~\ref{mai2} given by derivatives is
less meaningful. Nevertheless, we state it in this form since it
is directly related to~\eqref{bi3}.

The rest of the paper is organized as follows: in
Section~\ref{free} we prove Theorem~\ref{main2}. In Section~\ref{weight} we prove Theorem~\ref{main3}, which is essential for our bootstrap argument
in the proof of the main Theorem~\ref{main1} and also leads to the
rapid decay of eigenfunctions with nonzero eigenvalues
in~\eqref{bi22}. In Section~\ref{section5} we transfer the
limiting absorption principle~\eqref{bi3} from the case $L=0$
(Theorem~\ref{main2}) to the general case of admissible
perturbations, by means of the resolvent identity and Fredholm's
alternative. In Section~\ref{proof} we prove Theorem~\ref{main1}
and Corollary~\ref{mai2}. Finally, in Section~\ref{examples} we
prove Proposition~\ref{prop16}.

The authors would like to thank C. E. Kenig and B. Simon for useful discussions.

\section{Proof of Theorem \ref{main2}}\label{free}

For $\lambda\in[-\delta^{-1},-\delta]$ or $|\epsilon|>\delta$, we have
the elliptic bound
\begin{equation*}
||R_0(\lambda+i\epsilon)||_{W^{-1,2}\to W^{1,2}}\leq C_\delta,
\end{equation*}
which is stronger than \eqref{bu1}, in view of \eqref{tb2}. Thus,
we may assume that $\lambda\in[\delta,\delta^{-1}]$ and
$|\epsilon|\leq\delta$. Let
$\chi:\mathbb{R}^d\to[0,1]$ denote a smooth function supported in
the set $\{|\xi|\in[\sqrt{\lambda}/2,3\sqrt{\lambda}/2]\}$ and equal to $1$ in the set
$\{|\xi|\in[3\sqrt{\lambda}/4,5\sqrt{\lambda}/4]\}$. Let $\chi(D)$ and $(1-\chi)(D)$ denote the
operators defined by the Fourier multipliers $\xi\to\chi(\xi)$ and
$\xi\to1-\chi(\xi)$. For the operator
$(1-\chi)(D)R_0(\lambda+i\epsilon)$ we have again the stronger elliptic
bound
\begin{equation*}
||(1-\chi)(D)R_0(\lambda+i\epsilon)||_{W^{-1,2}\to W^{1,2}}\leq C_\delta.
\end{equation*}
It remains to prove that
\begin{equation*}
||\chi(D)R_0(\lambda+i\epsilon)||_{X\to X^\ast}\leq C_\delta,
\end{equation*}
which is equivalent to
\begin{equation}\label{tb20}
\begin{split}
&||S_1\chi(D)R_0(\lambda+i\epsilon)S_1||_{B\to B^\ast}+||S_{\frac{1}{d+1}}\chi(D)R_0(\lambda+i\epsilon)S_{\frac{1}{d+1}}||_{L^{p_d}\to L^{p_d'}}\\
&+||S_{1}\chi(D)R_0(\lambda+i\epsilon)S_{\frac{1}{d+1}}||_{L^{p_d}\to B^\ast}+||S_{\frac{1}{d+1}}\chi(D)R_0(\lambda+i\epsilon)S_{1}||_{B\to L^{p_d'}}\leq C_\delta.
\end{split}
\end{equation}
The $B\to
B^\ast$ bound in \eqref{tb20} follows from
\cite[Theorem 14.2.2]{Ho}. Also, since
\begin{equation*}
\langle S_{1/(d+1)}\chi(D)R_0(\lambda+i\epsilon)S_{1}f,g\rangle=\langle
f,S_{1}\chi(D)R_0(\lambda-i\epsilon)S_{1/(d+1)}g\rangle
\end{equation*}
for any $f,g\in\mathcal{S}(\mathbb{R}^d)$, it follows that
\begin{equation*}
||S_{1/(d+1)}\chi(D)R_0(\lambda+i\epsilon)S_{1}||_{B\to L^{p_d'}}=||S_{1}\chi(D)R_0(\lambda-i\epsilon)S_{1/(d+1)}||_{L^{p_d}\to B^\ast}.
\end{equation*}
Thus it remains to prove the $L^{p_d}\to L^{p_d'}$ and the
$L^{p_d}\to B^\ast$ bounds in \eqref{tb20}. The $L^{p_d}\to B^\ast$
bound follows from the work of A.~Ruiz and L.~Vega,
see~\cite[Theorem 3.1]{RuVe}. For the $L^{p_d}\to L^{p_d'}$ bound,
we notice that $S_{2/(d+1)}\chi(D)$ is a bounded operator on
$L^{p_d}$, so it suffices to prove that
\begin{equation*}
||R_0(\lambda+i\epsilon)||_{L^{p_d}\to L^{p'_d}}\leq C_\delta
\end{equation*}
uniformly in $\epsilon$ and $\lambda\in[\delta,\delta^{-1}]$. This follows from
\cite[Theorem 2.3]{KeRuSo}.

\section{Proof of Theorem \ref{main3}}\label{weight}

The constants $C_N$ in this section may depend on $N$ and the
dimension $d$, but not on $\gamma$.  We start with a lemma concerning the weight $\mu_{N,\gamma}$.

\newtheorem{lemma6}{Lemma}[section]
\begin{lemma6}\label{lemma6}
(a) We have
\begin{equation*}
\begin{split}
&\partial_{x_j}\mu_{N,\gamma}=\mu_{N,\gamma} b_j,\,\,\,\,\Delta\mu_{N,\gamma}=\mu_{N,\gamma} b;\\
&\partial_{x_j}\mu^{-1}_{N,\gamma}=-\mu^{-1}_{N,\gamma} b_j,\,\,\,\,\Delta\mu^{-1}_{N,\gamma}=\mu^{-1}_{N,\gamma}\widetilde{b},
\end{split}
\end{equation*}
for some functions $b_j$, $b$, and $\widetilde{b}$ with the
property that for any $x\in\mathbb{R}^n$
\begin{equation}\label{tb651}
\sum_{j=1}^d|b_j(x)|(1+|x|^2)^{1/2}+|b(x)|(1+|x|^2)+|\widetilde{b}(x)|(1+|x|^2)\leq C_N.
\end{equation}

(b) We have
\begin{equation*}
\frac{\mu_{N,\gamma}(x)}{\mu_{N,\gamma}(y)}+\frac{\mu_{N,\gamma}(y)}{\mu_{N,\gamma}(x)}\leq C_N(1+|x-y|^2)^N,
\end{equation*}
for any $x,y\in\mathbb{R}^d$.

(c) For any $r\in(0,\infty)$ and $x\in\mathbb{R}^d$,
\begin{equation*}
\mu_{N,\gamma}(x)\leq C_{N,r}\mu_{N,\gamma}(r x)\leq C_{N,r}\mu_{N,\gamma}(x).
\end{equation*}
\end{lemma6}

\begin{proof}[Proof of Lemma~\ref{lemma6}] The proof follows easily from the formula \eqref{tb3}.
\end{proof}

We also need a technical lemma that allows us to commute the
operators $S_\alpha$ and multiplication by the weight
$\mu_{N,\gamma}$.

\newtheorem{lemma6.1}[lemma6]{Lemma}
\begin{lemma6.1}\label{lemma6.1}
For $\alpha\in[-2,2]$
\begin{equation}\label{tb30}
||\mu S_\alpha\mu^{-1} S_{-\alpha}||_{L^p\to L^p}+||\mu S_\alpha\mu^{-1} S_{-\alpha}||_{B\to B}+||\mu S_\alpha\mu^{-1} S_{-\alpha}||_{B^\ast\to B^\ast}\leq C_N,
\end{equation}
and
\begin{equation}\label{tb31}
||S_\alpha\mu S_{-\alpha}\mu^{-1}||_{L^p\to L^p}+||S_\alpha\mu S_{-\alpha}\mu^{-1}||_{B\to B}+||S_\alpha\mu S_{-\alpha}\mu^{-1}||_{B^\ast\to B^\ast}\leq C_N,
\end{equation}
where $p\in\{p_d,2,p_d'\}$ and $\mu\in\{\mu_{N,\gamma},\mu_{N,\gamma}^{-1}\}$.
\end{lemma6.1}

\begin{proof}[Proof of Lemma~\ref{lemma6.1}] (a) By analytic interpolation,
it suffices to prove \eqref{tb30} and \eqref{tb31} for $\alpha=\pm2+i\beta$,
$\beta\in\mathbb{R}$, with constant $C_Ne^{\beta^2}$. Notice also that
\begin{equation*}
\langle\mu S_{-2-i\beta}\mu^{-1} S_{2+i\beta}f,\psi\rangle=\langle f,S_{2-i\beta}\mu^{-1}S_{-2+i\beta}\mu\psi\rangle,
\end{equation*}
for any $f\in\mathcal{S'}(\mathbb{R}^d)$ and
$\psi\in\mathcal{S}(\mathbb{R}^d)$. So it suffices to prove
\eqref{tb30} and \eqref{tb31} for $\alpha=2+i\beta$,
$\beta\in\mathbb{R}$. We use the fact that $S_2=-\Delta+1$ and
Lemma~\ref{lemma6}(a). Then, with $a_j=b_j$ and $a=b$, or
$a_j=-b_j$ and $a=\widetilde{b}$
\begin{equation}\label{tb40}
\mu S_{2+i\beta}\mu^{-1}S_{-2-i\beta}=\mu S_{i\beta}\mu^{-1}[S_{-i\beta}+2\sum_ja_j\partial_{x_j}S_{-2-i\beta}+aS_{-2-i\beta}]
\end{equation}
and
\begin{equation}\label{tb41}
S_{2+i\beta}\mu S_{-2-i\beta}\mu^{-1}=S_{i\beta}[\mu S_{-i\beta}\mu^{-1}+2\sum_ja_j\mu\partial_{x_j}S_{-2-i\beta}\mu^{-1}+a\mu S_{-2-i\beta}\mu^{-1}].
\end{equation}
The operator $S_{i\beta}$ is bounded on $B$, $B^\ast$, and $L^p$,
by the theory of singular integrals (for boundedness on $B$ and
$B^\ast$, notice that the kernel of this operator is rapidly
decreasing at $\infty$). The same is true for the operator
$2\sum_ja_j\partial_{x_j}S_{-2-i\beta}+aS_{-2-i\beta}$ in the
right-hand side of \eqref{tb40}, in view of Lemma~\ref{lemma6}(a).
Therefore it suffices to prove that if $m\in
C^\infty(\mathbb{R}^d)$ satisfies the differential bounds
\begin{equation}\label{tb43}
|\partial_\xi^\nu m(\xi)|\leq C_\nu(1+|\xi|^2)^{-|\nu|/2}
\end{equation}
for any $\xi\in\mathbb{R}^d$ and multi-index $\nu$, then
\begin{equation}\label{tb44}
||\mu m(D)\mu^{-1}||_{L^p\to L^p}+||\mu m(D)\mu^{-1}||_{B\to B}+||\mu m(D)\mu^{-1}||_{B^\ast\to B^\ast}\leq C_{N,m},
\end{equation}
where $p$ and $\mu$ are as in Lemma~\ref{lemma6.1}, and $m(D)$
denotes the operator defined by the Fourier multiplier $\xi\to
m(\xi)$.

To prove \eqref{tb44}, we use the fact the kernel of the operator
$m(D)$ has rapid decay away from the diagonal:
\begin{equation}\label{tb46}
|m(D)(x,y)|\leq C_\nu|x-y|^{-\nu}
\end{equation}
for any $x,y\in\mathbb{R}^d$ and integer $\nu\geq 0$. This follows
from \eqref{tb43} by integration by parts. Let $D_j$ denote the
sets in the definition of the spaces $B$ and $B^\ast$. We show
first that if $p\in\{p_d,2,p_d'\}$ then
\begin{equation}\label{tb45}
||\mathbf{1}_{D_{j'}}\mu m(D)\mu^{-1}\mathbf{1}_{D_j}||_{L^p\to L^p}\leq C_{N,m}2^{-|j-j'|},
\end{equation}
for any integers $j,j'\geq 0$. For $|j-j'|\geq 2$ we can simply
use \eqref{tb46} and the fact that the absolute value of the
kernel of the operator $\mathbf{1}_{D_{j'}}\mu
m(D)\mu^{-1}\mathbf{1}_{D_j}$ is $\mathbf{1}_{D_{j'}}(x)\mu(x)
|m(D)(x,y)|\mu^{-1}(y)\mathbf{1}_{D_j}(y)\leq
C_{N,\nu}(1+|x-y|)^{-\nu}2^{-|j-j'|}$, in view of Lemma~\ref{lemma6}(b).
For $|j-j'|\leq 1$ we use the fact that $m(D)$
defines a bounded operator on $L^p$:
\begin{equation*}
\begin{split}
||\mathbf{1}_{D_{j'}}\mu m(D)\mu^{-1}\mathbf{1}_{D_j}f||_{L^p}&\leq\sup_{x\in D_{j'}}\mu(x)||m(D)\mu^{-1}\mathbf{1}_{D_j}f||_{L^p}\\
&\leq\sup_{x\in D_{j'}}\mu(x)\sup_{y\in D_{j}}\mu^{-1}(y)||f||_{L^p},
\end{split}
\end{equation*}
which proves the $L^p\to L^p$ bound in \eqref{tb45}, in view of
Lemma~\ref{lemma6}(c).

We complete now the proof of \eqref{tb44}. For the $L^p\to L^p$ bound:
\begin{equation*}
\begin{split}
&||\mu m(D)\mu^{-1}f||_{L^p}^p=\sum_{j'=1}^\infty\big|\big|\sum_{j=1}^\infty\mathbf{1}_{D_{j'}}\mu m(D)\mu^{-1}\mathbf{1}_{D_j}f\big|\big|_{L^p}^p\\
&\leq C_{N,m}^p\sum_{j'=1}^\infty\big[\sum_{j=1}^\infty2^{-|j-j'|}||\mathbf{1}_{D_j}f||_{L^p}\big]^p\leq C_{N,m}^p\sum_{j=1}^\infty||\mathbf{1}_{D_j}f||_{L^p}^p=C_{N,m}^p||f||_{L^p}^p,
\end{split}
\end{equation*}
as desired. The proof of the $B\to B$ and $B^\ast\to B^\ast$ bounds
is similar, using the $L^2\to L^2$ bound in \eqref{tb45}. This
completes the proof of Lemma~\ref{lemma6.1}.
\end{proof}

For later use, we show that if $\chi\in C_0^\infty(\mathbb{R}^d)$
and $\chi(D)$ denotes the operator defined by the Fourier
multiplier $\xi\to\chi(\xi)$ then
\begin{equation}\label{tb60}
||\mu_{N,\gamma}\chi(D)g||_{X}\leq C_{N,\chi}||\mu_{N,\gamma}g||_{X}
\end{equation}
for any $g\in X$. For this, we notice first that if $g\in X$ then $\mu_{N,\gamma}g\in X$: let $g=g_1+g_2$ with $S_{-1/(d+1)}g_1\in L^{p_d}$ and $S_{-1}g_2\in B$. Then
\begin{equation*}
\begin{split}
||\mu_{N,\gamma}g||_X&\leq ||S_{-1/(d+1)}\mu_{N,\gamma}g_1||_{L^{p_d}}+||S_{-1}\mu_{N,\gamma}g_2||_{B}\\
&=||[S_{-1/(d+1)}\mu_{N,\gamma}S_{1/(d+1)}\mu_{N,\gamma}^{-1}]\mu_{N,\gamma}S_{-1/(d+1)}g_1||_{L^{p_d}}\\
&+||[S_{-1}\mu_{N,\gamma}S_{1}\mu_{N,\gamma}^{-1}]\mu_{N,\gamma}S_{-1}g_2||_{B}\\
&\leq C_{N,\gamma}[||S_{-1/(d+1)}g_1||_{L^{p_d}}+||S_{-1}g_2||_{B}]\leq C_{N,\gamma}||g||_X,
\end{split}
\end{equation*}
using Lemma~\ref{lemma6.1} and the fact that
$|\mu_{N,\gamma}(x)|\leq C_{N,\gamma}$. To prove \eqref{tb60}, let
$g=g_1+g_2$ be such  that $2||\mu_{N,\gamma}g||_{X}\geq
||S_{-1/(d+1)}\mu_{N,\gamma}g_1||_{L^{p_d}}+||S_{-1}\mu_{N,\gamma}g_2||_{B}$.
Then
\begin{equation*}
\begin{split}
&||\mu_{N,\gamma}\chi(D)g||_{X}\leq ||S_{-1/(d+1)}\mu_{N,\gamma}\chi(D)g_1||_{L^{p_d}}+||S_{-1}\mu_{N,\gamma}\chi(D)g_2||_{B}\\
&=||[S_{\frac{-1}{d+1}}\mu_{N,\gamma}S_{\frac{1}{d+1}}\mu_{N,\gamma}^{-1}][\mu_{N,\gamma}\chi(D)\mu_{N,\gamma}^{-1}][\mu_{N,\gamma}S_{\frac{-1}{d+1}}\mu_{N,\gamma}^{-1}S_{\frac{1}{d+1}}]S_{\frac{-1}{d+1}}\mu_{N,\gamma}g_1||_{L^{p_d}}\\
&+||[S_{-1}\mu_{N,\gamma}S_{1}\mu_{N,\gamma}^{-1}][\mu_{N,\gamma}\chi(D)\mu_{N,\gamma}^{-1}][\mu_{N,\gamma}S_{-1}\mu_{N,\gamma}^{-1}S_{1}]S_{-1}\mu_{N,\gamma}g_2||_{B}\\
&\leq C_{N,\chi}[||S_{-1/(d+1)}\mu_{N,\gamma}g_1||_{L^{p_d}}+||S_{-1}\mu_{N,\gamma}g_2||_{B}]\leq C_{N,\chi}||\mu_{N,\gamma}g||_{X},
\end{split}
\end{equation*}
using Lemma~\ref{lemma6.1} and \eqref{tb44}. This completes the proof of \eqref{tb60}.
\medskip

We turn now to the proof of Theorem~\ref{main3}. We prove first
the bound \eqref{bo1} under the additional restriction
\begin{equation}\label{tb50}
u\in\mathcal{S}(\mathbb{R}^d).
\end{equation}
In the case $\lambda\in[-\delta^{-1},-\delta]$ we prove the stronger elliptic bound
\begin{equation}\label{tb51}
||\mu_{N,\gamma}u||_{W^{1,2}}\leq C_{N,\delta}||\mu_{N,\gamma}(\Delta+\lambda)u||_{W^{-1,2}}.
\end{equation}
Let $f=S_{-1}\mu_{N,\gamma}(\Delta+\lambda)u$. Then
$S_1\mu_{N,\gamma}u=S_1\mu_{N,\gamma}(\Delta+\lambda)^{-1}\mu_{N,\gamma}^{-1}S_1f$.
Thus
\begin{equation*}
\begin{split}
&||\mu_{N,\gamma}u||_{W^{1,2}}=||S_1\mu_{N,\gamma}u||_{L^2}\\
&=||[S_1\mu_{N,\gamma}S_{-1}\mu_{N,\gamma}^{-1}][\mu_{N,\gamma}S_1(\Delta+\lambda)^{-1}S_1\mu_{N,\gamma}^{-1}][\mu_{N,\gamma}S_{-1}\mu_{N,\gamma}^{-1}S_1]f||_{L^2}\\
&\leq C_N||f||_{L^2},
\end{split}
\end{equation*}
using Lemma~\ref{lemma6.1} and \eqref{tb44}. This proves \eqref{tb51}.

The proof in the case $\lambda\in[\delta,\delta^{-1}]$ is more
difficult, since the elliptic bound \eqref{tb51} does not hold.
For some small constant $\varepsilon_0=\varepsilon_0(\delta)>0$
(to be fixed later), let $A=\{\xi^1,\ldots,\xi^m\}$ denote a
$\varepsilon_0/100$-net on the sphere $\{|\xi|=\sqrt\lambda\}$ in
the Fourier space. Using this net, we construct a partition of $1$
in the Fourier space. We have $1=\chi_0+\chi_1+\ldots+\chi_m$,
where $ \chi_0,\chi_1,\ldots,\chi_m:\mathbb{R}^d\to[0,1]$ are
smooth functions, $\chi_0$ is supported in the set
$\{\xi:|\sqrt\lambda-|\xi|\,|\geq \varepsilon_0/10\}$, and
$\chi_j$ is supported in the set $\{\xi:|\xi^j-\xi|\,\leq
\varepsilon_0/2\}$. Let $\chi_j(D)$ denote the operator defined by
the Fourier multiplier $\chi_j$.

An estimate similar to the proof of \eqref{tb51} shows that
\begin{align}
||\mu_{N,\gamma}\chi_0(D)u||_{X^\ast}&\leq C||\mu_{N,\gamma}\chi_0(D)u||_{W^{1,2}}\leq C_{N,\delta}||\mu_{N,\gamma}\chi_0(D)(\Delta+\lambda)u||_{W^{-1,2}}\nonumber \\
&\leq C_{N,\delta}||\mu_{N,\gamma}(\Delta+\lambda)u||_{X},\label{tb62}
\end{align}
using \eqref{tb60}. It remains to prove a similar estimate for
$||\mu_{N,\gamma}\chi_j(D)u||_{X^\ast}$, $j=1,\ldots,m$. Let
$\widetilde{\chi}_j$ denote a smooth function supported in the set
$\{\xi:|\xi^j-\xi|\,\leq \varepsilon_0\}$ and equal to $1$ in the
set $\{\xi:|\xi^j-\xi|\,\leq \varepsilon_0/2\}$. Thus
$\chi_j\widetilde{\chi}_j=\chi_j$. Then
\begin{equation}\label{tb65}
\begin{split}
&||\mu_{N,\gamma}\chi_j(D)u||_{X^\ast}\leq ||S_{1/(d+1)}\mu_{N,\gamma}\chi_j(D)u||_{L^{p'_d}}+||S_{1}\mu_{N,\gamma}\chi_j(D)u||_{B^\ast}\\
&=||[S_{1/(d+1)}\mu_{N,\gamma}S_{-1/(d+1)}\mu_{N,\gamma}^{-1}][\mu_{N,\gamma}S_{1/(d+1)}\widetilde{\chi}_j(D)\mu_{N,\gamma}^{-1}]\mu_{N,\gamma}\chi_j(D)u||_{L^{p'_d}}\\
&+||[S_{1}\mu_{N,\gamma}S_{-1}\mu_{N,\gamma}^{-1}][\mu_{N,\gamma}S_{1}\widetilde{\chi}_j(D)\mu_{N,\gamma}^{-1}]\mu_{N,\gamma}\chi_j(D)u||_{B^\ast}\\
&\leq C_{N,\delta}||\mu_{N,\gamma}\chi_j(D)u||_{L^{p'_d}\cap B^\ast},
\end{split}
\end{equation}
using Lemma~\ref{lemma6.1} and~\eqref{tb44}.
A similar estimate,
using again Lemma~\ref{lemma6.1} and~\eqref{tb44}, together with
decompositions as in the proof of \eqref{tb60}, shows that
\begin{equation}\label{tb66}
||\mu_{N,\gamma}\chi_j(D)(\Delta+\lambda)u||_{L^{p_d}+B}\leq C_{N,\delta}||\mu_{N,\gamma}\chi_j(D)(\Delta+\lambda)u||_{X}.
\end{equation}
The estimates \eqref{tb60}, \eqref{tb62}, \eqref{tb65}, and \eqref{tb66} show that it suffices to prove that for any $u\in\mathcal{S}(\mathbb{R}^d)$
\begin{equation}\label{bo10}
||\mu_{N,\gamma}\chi_j(D)u||_{L^{p'_d}\cap B^\ast}\leq
C_{N,\delta}||\mu_{N,\gamma}\chi_j(D)(\Delta+\lambda)u||_{L^{p_d}+B},j=1,\ldots,m.
\end{equation}

It remains to prove \eqref{bo10}. By rescaling and rotation
invariance, we may assume that $\lambda=1$ and $\chi_j=\chi$ is a smooth function
supported in the ball of radius $\varepsilon_0/2$ around the unit
vector $\xi^+=(0,\ldots,0,1)$. Let $\xi^+_1,\ldots,\xi^+_d$ denote a basis of
$\mathbb{R}^d$ of unit vectors in the ball
$\{\xi:|\xi-\xi^+|\leq\varepsilon_0/2\}$. Clearly
\begin{equation*}
|x|\leq C(|x\cdot\xi_1^+|+\ldots+|x\cdot\xi_d^+|).
\end{equation*}
It follows easily that
\begin{equation*}
\mu_{N,\gamma}(x)\approx[\widetilde{\mu}_{N,\gamma}(|x\cdot\xi_1^+|)+\ldots+\widetilde{\mu}_{N,\gamma}(|x\cdot\xi_d^+|)],
\end{equation*}
where, for $t\in[0,\infty)$
\begin{equation*}
\widetilde{\mu}_{N,\gamma}(t)=\frac{(1+t^2)^N}{(1+\gamma t^2)^N}.
\end{equation*}
Thus we may replace the weight $\mu_{N,\gamma}(x)$ in
\eqref{bo10} with $\widetilde{\mu}_{N,\gamma}(|x\cdot\xi_l^+|)$ (in both
sides of the inequality). To summarize, by rotation invariance, it
remains to prove that
\begin{equation}\label{bo20}
||\widetilde{\mu}_{N,\gamma}(|x_d|)u||_{L^{p'_d}\cap B^\ast}\leq
C_{N}||\widetilde{\mu}_{N,\gamma}(|x_d|)(\Delta+1)u||_{L^{p_d}+B},
\end{equation}
for all functions $u\in\mathcal{S}(\mathbb{R}^d)$ with the
property that $\hat{u}$ is supported in the ball
$\{\xi:|\xi-\xi^+|\leq\varepsilon_0\}$.

Let $\mathbf{1}_+$ and $\mathbf{1}_-$
denote the characteristic functions of the intervals $[0,\infty)$
and $(-\infty,0)$. Let $f=(\Delta+1)u$,
$F(x)=\widetilde{\mu}_{N,\gamma}(|x_d|)f(x)$, and
$U(x)=\widetilde{\mu}_{N,\gamma}(|x_d|)u(x)$. Let $\widetilde{u}(\xi',x_d)$
and $\widetilde{f}(\xi',x_d)$ denote the partial Fourier
transforms of the functions $u$ and $f$ in the variable
$x'=(x_1,\ldots,x_{d-1})$. The equation $f=(\Delta+1)u$ is
equivalent to
\begin{equation}\label{bo21}
[\partial_{x_d}^2+(1-|\xi'|^2)]\widetilde{u}(\xi',x_d)=\widetilde{f}(\xi',x_d).
\end{equation}
The functions $\widetilde{u}$ and $\widetilde{f}$ are supported in
the ball $\{\xi':|\xi'|\leq\varepsilon_0\ll 1\}$. By integration by parts,
\begin{equation}\label{bo40}
\begin{split}
\widetilde{u}(\xi',x_d)&=-\int_{x_d}^\infty \widetilde{f}(\xi',y_d)\frac{\sin(\sqrt{1-|\xi'|^2}(x_d-y_d))}{\sqrt{1-|\xi'|^2}}\,dy_d\\
&=\int_{-\infty}^{x_d} \widetilde{f}(\xi',y_d)\frac{\sin(\sqrt{1-|\xi'|^2}(x_d-y_d))}{\sqrt{1-|\xi'|^2}}\,dy_d.
\end{split}
\end{equation}
We use the formula in the first line of \eqref{bo40} when $x_d\geq 0$, and the formula in the second line when $x_d\leq 0$. Let $\phi:\mathbb{R}^{d-1}\to[0,1]$ denote a smooth function supported in the ball $\{\xi':|\xi'|\leq 2\varepsilon_0\}$ and equal to $1$ in the ball $\{\xi':|\xi'|\leq\varepsilon_0\}$. By taking the inverse Fourier transform in the variable $\xi'$ we have
\begin{equation}\label{bo42}
\begin{split}
u(x',x_d)&=c\mathbf{1}_+(x_d)\int_{\mathbb{R}^d}f(y',y_d)\mathbf{1}_-(x_d-y_d)H(x'-y',x_d-y_d)\,dy\\
&-c\mathbf{1}_-(x_d)\int_{\mathbb{R}^d}f(y',y_d)\mathbf{1}_+(x_d-y_d)H(x'-y',x_d-y_d)\,dy,
\end{split}
\end{equation}
where
\begin{equation}\label{bo41}
H(z',z_d)=\int_{\mathbb{R}^{d-1}}e^{iz'\cdot\xi'}\frac{\sin(\sqrt{1-|\xi'|^2}z_d)}{\sqrt{1-|\xi'|^2}}\phi(\xi')\,d\xi'.
\end{equation}
By multiplying with the weight $\widetilde{\mu}_{N,\gamma}$, we have
\begin{equation}\label{bo45}
U(x',x_d)=c\int_{\mathbb{R}^d}F(y',y_d)K(x',x_d,y',y_d)\,dy'dy_d,
\end{equation}
where
\begin{equation}\label{bo46}
\begin{split}
K(x',x_d,y',y_d)=[\mathbf{1}_+(x_d)&\mathbf{1}_-(x_d-y_d)-\mathbf{1}_-(x_d)\mathbf{1}_+(x_d-y_d)]\\
&\frac{\widetilde{\mu}_{N,\gamma}(|x_d|)}{\widetilde{\mu}_{N,\gamma}(|y_d|)}H(x'-y',x_d-y_d).
\end{split}
\end{equation}
It is important to notice that $K(x',x_d,y',y_d)=0$ if $|y_d|>|x_d|$; therefore the weight $\widetilde{\mu}_{N,\gamma}(|x_d|)/\widetilde{\mu}_{N,\gamma}(|y_d|)$ is always $\leq 1$. Let $T$ denote the operator defined by the kernel $K$ in the right-hand side of \eqref{bo45}. It remains to prove that $T$ extends to a bounded operator from $L^{p_d}+B$ to $L^{p'_d}\cap B^\ast$.

\newtheorem{prop1}[lemma6]{Lemma}
\begin{prop1}\label{prop1}
For $f\in\mathcal{S}(\mathbb{R}^d)$ we have
\begin{equation*}
||Tf||_{B^\ast}\leq C||f||_{B}.
\end{equation*}
\end{prop1}
\begin{proof}[Proof of Lemma~\ref{prop1}] This is essentially proved in \cite[Chapter XIV]{Ho}.
We need the observation that
\begin{equation*}
||Tf||_{B^\ast}\leq C\sup_{x_d\in\mathbb{R}}||Tf(.,x_d)||_{L^2_{x'}}\text{ and }||f||_{B}\geq C^{-1}\int_{\mathbb{R}}||f(.,y_d)||_{L^2_{y'}}\,dy_d.
\end{equation*}
In addition, if $x_d\geq 0$ then
\begin{equation*}
\begin{split}
||Tf(.,x_d)||_{L^2_{x'}}&=||\int_{\mathbb{R}^d}f(y',y_d)\mathbf{1}_-(x_d-y_d)\frac{\widetilde{\mu}_{N,\gamma}(x_d)}{\widetilde{\mu}_{N,\gamma}(y_d)}H(x'-y',x_d-y_d)\,dy||_{L^2_{x'}}\\
&\leq\int_{x_d}^\infty\frac{\widetilde{\mu}_{N,\gamma}(x_d)}{\widetilde{\mu}_{N,\gamma}(y_d)}||\int_{\mathbb{R}^{d-1}}f(y',y_d)H(x'-y',x_d-y_d)\,dy'||_{L^2_{x'}}dy_d\\
&\leq C\int_{\mathbb{R}}||f(.,y_d)||_{L^2_{y'}}\,dy_d,
\end{split}
\end{equation*}
where the last inequality follows from Plancherel's Theorem and the monotonicity of the weight $\widetilde{\mu}_{N,\gamma}$. The estimate in the case $x_d<0$ is similar.
\end{proof}

\newtheorem{prop2}[lemma6]{Lemma}
\begin{prop2}\label{prop2}
For $f\in\mathcal{S}(\mathbb{R}^d)$ we have
\begin{equation}\label{bo60}
||Tf||_{L^{p_d'}}\leq C||f||_{L^{p_d}}.
\end{equation}
\end{prop2}
\begin{proof}[Proof of Lemma~\ref{prop2}] Let
\begin{equation*}
W_N(x_d,y_d)=[\mathbf{1}_+(x_d)\mathbf{1}_-(x_d-y_d)-\mathbf{1}_-(x_d)\mathbf{1}_+(x_d-y_d)]\frac{\widetilde{\mu}_{N,\gamma}(|x_d|)}{\widetilde{\mu}_{N,\gamma}(|y_d|)}
\end{equation*}
denote the weight in the definition \eqref{bo46} of the kernel
$K$, $W_N(x_d,y_d)\in[-1,1]$. We use analytic interpolation. For
$\sigma\in\mathbb{C}$, $\Re\sigma\in[-(d-1)/2,1]$, let
\begin{equation*}
K^\sigma(x',x_d,y',y_d)=e^{\sigma^2}(1-\sigma)(1+|x_d-y_d|)^{-\sigma}K(x',x_d,y',y_d),
\end{equation*}
and $T^\sigma$ the operator defined by the kernel $K^\sigma$. By
analytic interpolation, it suffices to prove that
\begin{equation}\label{bo55}
||T^\sigma||_{L^1\to L^\infty}\leq C\,\text{ if }\Re\sigma=-(d-1)/2,
\end{equation}
and
\begin{equation}\label{bo56}
||T^\sigma||_{L^2\to L^2}\leq C\,\text{ if }\Re\sigma=1.
\end{equation}

The bound \eqref{bo55} follows easily since $|H(z',z_d)|\leq
C(1+|z_d|)^{-(d-1)/2}$, by stationary phase arguments.

To prove \eqref{bo56}, we take partial Fourier transforms in the
variables $y'$ and $x'$. Let
$W_N^\sigma(x_d,y_d)=e^{\sigma^2}(1-\sigma)(1+|x_d-y_d|)^{-\sigma}W_N(x_d,y_d)$.
An easy computation shows that
\begin{equation}\label{bo66}
\widetilde{T^\sigma f}(\eta',x_d)=c\int_{\mathbb{R}}\widetilde{f}(\eta',y_d)W_N^\sigma(x_d,y_d)\frac{\sin[\sqrt{1-|\eta'|^2}(x_d-y_d)]}{\sqrt{1-|\eta'|^2}}\phi(\eta')\,dy_d.
\end{equation}
Notice that
\begin{equation*}
\sin[\sqrt{1-|\eta'|^2}(x_d-y_d)]=c[e^{i\sqrt{1-|\eta'|^2}(x_d-y_d)}-e^{-i\sqrt{1-|\eta'|^2}(x_d-y_d)}].
\end{equation*}
We substitute this into \eqref{bo66}. Notice that the exponential
factors $e^{\pm i\sqrt{1-|\eta'|^2}x_d}$ and $e^{\pm
i\sqrt{1-|\eta'|^2}y_d}$ can be paired with $\widetilde{T^\sigma
f}(\eta',x_d)$ and $\widetilde{f}(\eta',y_d)$ respectively. By
Plancherel's theorem, the $L^2$ bound \eqref{bo56} would follow once
we prove that the kernel $W_N^\sigma$ defines a bounded operator
on $L^2(\mathbb{R})$:
\begin{equation}\label{bo67}
||\int_\mathbb{R}h(y_d)W_N^\sigma(x_d,y_d)\,dy_d||_{L^2_{x_d}}\leq C||h||_{L^2_{y_d}}\,\text{ if }\Re\sigma=1.
\end{equation}
We will use the maximal operator
\begin{equation*}
Mh(t)=\sup_{r\in\mathbb{R}}\left|\int_0^rh(t-s)e^{\sigma^2}(1-\sigma)(1+|s|)^{-\sigma}\,ds\right|.
\end{equation*}
For $h\in\mathcal{S}(\mathbb{R})$ and $x_d\geq 0$ we have
\begin{equation}\label{bo70}
\begin{split}
|\int_\mathbb{R}h(y_d)&W_N^\sigma(x_d,y_d)\,dy_d|
=\!|\!\!\int_{x_d}^\infty
h(y_d)e^{\sigma^2}(1-\sigma)(1+|x_d-y_d|)^{-\sigma}\frac{\widetilde{\mu}_{N,\gamma}(x_d)}{\widetilde{\mu}_{N,\gamma}(y_d)}\,dy_d|\\
&=|\int_{x_d}^\infty\frac{d}{dy_d}\left(\int_{x_d}^{y_d}h(s)e^{\sigma^2}(1-\sigma)(1+|x_d-s|)^{-\sigma}\,ds\right)\frac{\widetilde{\mu}_{N,\gamma}(x_d)}{\widetilde{\mu}_{N,\gamma}(y_d)}\,dy_d|\\
&=|\int_{x_d}^\infty\left(\int_{x_d}^{y_d}h(s)e^{\sigma^2}(1-\sigma)(1+|x_d-s|)^{-\sigma}\,ds\right)\frac{d}{dy_d}\frac{\widetilde{\mu}_{N,\gamma}(x_d)}{\widetilde{\mu}_{N,\gamma}(y_d)}\,dy_d|\\
&\leq\int_{x_d}^\infty
Mh(x_d)\left|\frac{d}{dy_d}\frac{\widetilde{\mu}_{N,\gamma}(x_d)}{\widetilde{\mu}_{N,\gamma}(y_d)}\right|\,dy_d\leq
Mh(x_d).
\end{split}
\end{equation}
The last inequality is due to the fact that the function
$y_d\to[\mu_{N,\gamma}(x_d)/\mu_{N,\gamma}(y_d)]$ is
nonincreasing, thus it has bounded variation. A similar
computation proves the estimate \eqref{bo70} in the case $x_d<0$.
In addition, when $\Re\sigma=1$, the kernels
$\chi_{\pm}(s)e^{\sigma^2}(1-\sigma)(1+|s|)^{-\sigma}$ are
Calder\'{o}n--Zygmund kernels, uniformly in $\sigma$. Therefore
the maximal operator $M$ is bounded on $L^2(\mathbb{R})$ (see, for
example, \cite[Chapter I, Section 7]{St}, so the bound
\eqref{bo67} follows from \eqref{bo70}.
\end{proof}

\newtheorem{prop3}[lemma6]{Lemma}
\begin{prop3}\label{prop3}
For $f\in\mathcal{S}(\mathbb{R}^d)$ we
have
\begin{equation}\label{bo80}
||T||_{L^{p_d}\to B^\ast}+||T||_{B\to L^{p'_d}}\leq C.
\end{equation}
\end{prop3}

\begin{proof}[Proof of Lemma~\ref{prop3}] We prove the bound for the first term in \eqref{bo80} (the proof for the second term is identical). As in the proof of Lemma \ref{prop1}, it suffices to prove that
\begin{equation}\label{pp1}
||Tf(.,x_d)||_{L^2_{x'}}\leq C||f||_{L^{p_d}}.
\end{equation}
Assuming $x_d$ fixed, let $g(y',y_d)=f(y',y_d)W_N(x_d,y_d)$, where $W_N(x_d,y_d)\in[-1,1]$ is the weight defined in the proof of Lemma \ref{prop2}. Clearly, $||g||_{L^{p_d}}\leq ||f||_{L^{p_d}}$. Also,
\begin{equation*}
\begin{split}
Tf(x',&x_d)=\int_{\mathbb{R}^d}g(y',y_d)H(x'-y',x_d-y_d)\,dy'dy_d=\frac{1}{2i}\int_{\mathbb{R}^{d-1}}e^{ix'\cdot\xi'}\frac{\phi(\xi')}{\sqrt{1-|\xi'|^2}}\\
&[\widehat{g}(\xi',\sqrt{1-|\xi'|^2})e^{i\sqrt{1-|\xi'|^2}x_d}-\widehat{g}(\xi',-\sqrt{1-|\xi'|^2})e^{-i\sqrt{1-|\xi'|^2}x_d}]\,d\xi'.
\end{split}
\end{equation*}
The bound \eqref{pp1} follows from Plancherel theorem and the Stein-Tomas restriction theorem.
\end{proof}
\medskip

We remove now the restriction \eqref{tb50}. Let
$\varphi:\mathbb{R}^d\to[0,C]$ denote a smooth function supported
in the ball $\{x:|x|\leq 1\}$ with
$\int_{\mathbb{R}^d}\varphi\,dx=1$, and
$\chi:\mathbb{R}^d\to[0,1]$ a smooth function supported in the
ball $\{x:|x|\leq 2\}$ and equal to $1$ in the ball $\{x:|x|\leq
1\}$. For $\varepsilon\in(0,1]$ and $r\in[1,\infty)$, let
$\varphi_\varepsilon(x)=\varepsilon^{-d}\varphi(x/\varepsilon)$
and $\chi_r(x)=\chi(x/r)$. Let
\begin{equation*}
u_{\varepsilon,r}(x)=\chi_r(x)(u\ast\varphi_\varepsilon)(x).
\end{equation*}
Clearly $u_{\varepsilon,r}\in\mathcal{S}(\mathbb{R}^d)$ and
\begin{equation*}
(\Delta+\lambda)u_{\varepsilon,r}=\chi_r[(\Delta+\lambda)u\ast\varphi_\varepsilon]+2\nabla\chi_r\cdot\nabla(u\ast\varphi_\varepsilon)+\Delta\chi_r(u\ast\varphi_\varepsilon).
\end{equation*}
We apply Theorem~\ref{main3} to the Schwartz function $u_{\varepsilon,r}$. The result is
\begin{equation}\label{bo101}
\begin{split}
||\mu_{N,\gamma}&\chi_r(u\ast\varphi_\varepsilon)||_{X^\ast}\leq C_{N,\delta}||\mu_{N,\gamma}\chi_r[(\Delta+\lambda)u\ast\varphi_\varepsilon]||_{X_q}\\
&+C_{N,\delta}[||\mu_{N,\gamma}\nabla\chi_r\cdot\nabla(u\ast\varphi_\varepsilon)||_B+||\mu_{N,\gamma}\Delta\chi_r(u\ast\varphi_\varepsilon)||_B].
\end{split}
\end{equation}
The function $\nabla\chi_r$ and $\Delta\chi_r$ are both supported
in the set $\{x:|x|\in[r,2r]\}$ and dominated by $C/r$. In
addition, $|\nabla(u\ast\varphi_\varepsilon)|\leq
C\varepsilon^{-1}
|u|\ast(\varepsilon^{-d}|\nabla\varphi(./\varepsilon)|)$. By
\eqref{bo100}, assuming $\varepsilon$ fixed and letting
$r\to\infty$ in \eqref{bo101}, we have
\begin{equation*}
||\mu_{N,\gamma}(u\ast\varphi_\varepsilon)||_{X^\ast}\leq
C_{N,\delta}||\mu_{N,\gamma}(\Delta+\lambda)v\ast\varphi_\varepsilon||_{X}.
\end{equation*}
The theorem follows by letting $\varepsilon\to 0$.

\section{The operators $\mathrm{Id}_{X^\ast}+R_0(\lambda\pm i\epsilon)L$}\label{section5}

Following the classical scheme, we will transfer some of the previous
estimates for the free resolvent to the perturbed resolvent by means
of the resolvent identity. This requires inverting
\begin{equation}
\label{eq:fred}
 {\mathrm Id}_{X^\ast} + R_0(\lambda\pm i0)L
\end{equation}
as an operator on $X^*$. We start with the definition of the
operators $R_0(\lambda\pm i0)$.

\newtheorem{lemma65}{Lemma}[section]
\begin{lemma65}\label{lemma65}(a) The map $z\to R_0(z)$ defines an
analytic map from $\mathbb{C}\setminus[0,\infty)$ to
$\mathcal{L}(X,X^\ast)$.

(b) For any $\lambda\in(0,\infty)$ there are operators $R_0(\lambda+i0),R_0(\lambda-i0)\in\mathcal{L}(X,X^\ast)$ with the property that
\begin{equation}\label{pp40}
||R_0(\lambda\pm i0)||_{X\to X^\ast}\leq C_\delta\text{ for any }\lambda\in[\delta,\delta^{-1}],\delta>0.
\end{equation}
In addition, for any sequences $\{\lambda_n\}_{n=1}^\infty\subseteq(0,\infty)$ and $\{\epsilon_n\}_{n=1}^\infty\subseteq[0,\infty)$ with $\lambda_n\to\lambda$ and $\epsilon_n\to0$, we have
\begin{equation}\label{pp41}
\lim_{n\to\infty}\langle R_0(\lambda_n\pm i\epsilon_n)f,\phi\rangle=\langle R_0(\lambda\pm i0)f,\phi\rangle
\end{equation}
for any $f\in X$ and $\phi\in\mathcal{S}(\mathbb{R}^d)$, and
\begin{equation}\label{pp42}
\lim_{n\to\infty}||\mathbf{1}_{\{|x|\leq R\}}[R_0(\lambda_n\pm i\epsilon_n)f-R_0(\lambda\pm i0)f]||_{L^2}=0
\end{equation}
for any $f\in X$ and $R\geq 1$.

(c) For $\lambda\in\mathbb{R}\setminus\{0\}$, $\epsilon\geq 0$,
and $g\in X$
\begin{equation*}
[-\Delta-(\lambda\pm i\epsilon)]R_0(\lambda\pm i\epsilon)g=g
\end{equation*}
in the sense of distributions (by a slight abuse of notation we
let $R_0(\lambda+i0)=R_0(\lambda-i0):=R_0(\lambda)$ when
$\lambda\in(-\infty,0)$).
\end{lemma65}
\begin{proof}[Proof of Lemma~\ref{lemma65}] Part (a) follows directly
from the definitions; in fact, the map $z\to R_0(z)$ defines an
analytic map from $\mathbb{C}\setminus[0,\infty)$ to
$\mathcal{L}(W^{-1,2},W^{1,2})$.

For part (b), we use the fact that $R_0(z)f=f\ast R_z$ for $z\in\mathbb{C}\setminus[0,\infty)$ and $f\in\mathcal{S}(\mathbb{R}^d)$, where
\begin{equation}
\label{eq:res0}
R_z(x)=C(z^{1/2}/|x|)^{(d-2)/2}
K_{(d-2)/2} (-iz^{1/2} |x|).
\end{equation}
Here $K_\nu$ denote the Bessel potentials and, as before, $\Im(z^{1/2})>0$
(see \cite[p.~288]{GeSh}). Standard estimates on the Bessel potentials show that if $|z|\in[\delta,\delta^{-1}]$ then
\begin{equation}\label{pp45}
|R_z(x)|\leq C_\delta
\begin{cases}
|x|^{-(d-1)/2}\quad&\text{ if }|x|\geq 1;\\
|x|^{-(d-2)}\quad&\text{ if }|x|\leq 1\text{ and }d\geq 3;\\
\log(2/|x|)\quad&\text{ if }|x|\leq 1\text{ and }d=2.
\end{cases}
\end{equation}
We define the kernels $R_{\lambda+i0}(x)$ and $R_{\lambda-i0}(x)$ using the formula \eqref{eq:res0} and letting $z\to\lambda+i0$ and $z\to\lambda-i0$. The kernels $R_{\lambda\pm i0}(x)$ satisfy the bound \eqref{pp45}. Then, for $f\in\mathcal{S}(\mathbb{R}^d)$, we define
\begin{equation*}
R_0(\lambda\pm i0)f:=f\ast R_{\lambda\pm i0}.
\end{equation*}
Using the Lebesgue dominated convergence theorem and \eqref{pp45},
\begin{equation}\label{pp49}
\lim_{n\to\infty}f\ast R_{\lambda_n\pm i\epsilon_n}(x)=f\ast R_{\lambda\pm i0}(x)
\end{equation}
for any $f\in\mathcal{S}(\mathbb{R}^d)$ and $x\in\mathbb{R}^d$, where $\lambda_n$ and $\epsilon_n$ are as in part (b). Using the Fatou lemma and Theorem \ref{main2}, for any $f\in\mathcal{S}(\mathbb{R}^d)$ and $\lambda\in[\delta,\delta^{-1}]$,
\begin{equation*}
\begin{split}
||R_0(\lambda\pm i0)f||_{S_{-1}(B^\ast)}&=||(S_1f)\ast R_{\lambda\pm i0}||_{B^\ast}\leq\limsup_{n\to\infty}||(S_1f)\ast R_{\lambda\pm i/n}||_{B^\ast}\\
&\leq \limsup_{n\to\infty}||f\ast R_{\lambda\pm i/n}||_{S_{-1}(B^\ast)}\leq C_{\delta}||f||_X.
\end{split}
\end{equation*}
A similar estimate shows that $||R_0(\lambda\pm i0)f||_{W^{1/(d+1),p_d'}}\leq C_{\delta}||f||_X$ for any $f\in\mathcal{S}(\mathbb{R}^d)$. Thus the operators $R_0(\lambda\pm i0)$ extend to bounded operators from $X$ to $X^\ast$ and \eqref{pp40} holds.

To prove the limits \eqref{pp41} and \eqref{pp42}, we notice that
we may assume $f\in\mathcal{S}(\mathbb{R}^d)$, in view of
\eqref{pp40} and the fact that $\mathcal{S}(\mathbb{R}^d)$ is
dense in $X$. The limits \eqref{pp41} and \eqref{pp42} then follow
from the Lebesgue dominated convergence theorem, the poinwise
limit \eqref{pp49}, and the observation that $|f\ast
R_{\lambda_n\pm i\epsilon_n}(x)|\leq
C||f||_{\mathcal{S}(\mathbb{R}^d)}$ (using \eqref{pp45}).

For part (c) we have to prove that for $g\in X$ and
$\phi\in\mathcal{S}(\mathbb{R}^d)$,
\begin{equation*}
\langle[-\Delta-(\lambda\pm i\epsilon)]R_0(\lambda\pm i\epsilon)g,\phi\rangle=\langle g,\phi\rangle,
\end{equation*}
which is equivalent to
\begin{equation}\label{new2}
\langle R_0(\lambda\pm i\epsilon)g,[-\Delta-(\lambda\mp i\epsilon)]\phi\rangle=\langle g,\phi\rangle.
\end{equation}
When $\lambda<0$ or $\epsilon\neq 0$, the identity \eqref{new2} is
clear, for any $g\in\mathcal{S}'(\mathbb{R}^d)$. When $\lambda>0$
and $\epsilon=0$, the identity \eqref{new2} follows from
\eqref{pp41}.
\end{proof}

Next we establish the compactness of the operator $L$.

\newtheorem{lemma7}[lemma65]{Lemma}
\begin{lemma7}\label{lemma7}
If $L$ is an admissible perturbation, then $L:X^\ast\to X$ is a
compact operator.
\end{lemma7}
\begin{proof} [Proof of Lemma~\ref{lemma7}] Let
$\{f_n\}_{n=1}^\infty$ denote a sequence of functions in $X^\ast$
with $||f_n||_{X^\ast}\leq 1$. Let $\chi:\mathbb{R}^d\to[0,1]$
denote a smooth function supported in the set $\{|x|\leq 2\}$ and
equal to $1$ in the set $\{|x|\le 1\}$. For any $r\geq 1$ let
$\chi_r(x)=\chi(x/r)$. Since $||S_1(f_n)||_{B^\ast}\leq 1$,
\begin{equation}\label{tb90}
||f_n\chi_r||_{W^{1,2}}\leq C_r
\end{equation}
for any $n\geq 1$ and $r\geq 1$.

We use first \eqref{tb90} with $r=1$. By the Rellich--Kondrachov
compactness theorem, there is a subsequence
$\{f_{1,n}\}_{n=1}^\infty\subseteq\{f_n\}_{n=1}^\infty$ and a
function $g_1\in L^2$ with the property that
\begin{equation*}
\lim_{n\to\infty}f_{1,n}\chi_1=g_1\text{ in }L^2.
\end{equation*}
We repeat this argument inductively for $r=2,3,\ldots$ and
construct subsequences
$\{f_{k,n}\}_{n=1}^\infty\subseteq\{f_{k-1,n}\}_{n=1}^\infty$ and
functions $g_k\in L^2$ with the property that
\begin{equation}\label{tb91}
\lim_{n\to\infty}f_{k,n}\chi_k=g_k\text{ in }L^2.
\end{equation}

We consider the diagonal subsequence $\widetilde{f}_k:=f_{k,k}$,
$k=1,2,\ldots$. It remains to prove that $L\widetilde{f}_k$ is a
Cauchy sequence in $X$. Given $\varepsilon>0$, we use \eqref{tb6}
with $N=0$. Therefore, there are constants $A_\varepsilon$ and
$R_\varepsilon$ with the property that
\begin{equation*}
||L(\widetilde{f}_k-\widetilde{f}_{k'})||_X\leq(\varepsilon/4)
||\widetilde{f}_k-\widetilde{f}_{k'}||_{X^\ast}+A_\varepsilon
||(\widetilde{f}_k-\widetilde{f}_{k'})\chi_{R\varepsilon}||_{L^2}.
\end{equation*}
By \eqref{tb91} and the definition of $\widetilde{f}_k$
\begin{equation*}
\limsup_{k,k'\to\infty}||(\widetilde{f}_k-\widetilde{f}_{k'})\chi_{R\varepsilon}||_{L^2}=0.
\end{equation*}
Thus $||L(\widetilde{f}_k-\widetilde{f}_{k'})||_X\leq\varepsilon$
for $k,k'$ large enough, as desired.
\end{proof}

The following is a technical lemma which will be needed in the
proof of invertibility of ${\rm Id_{X^*}}+R_0(\lambda\pm i\epsilon)L$.

\newtheorem{lemma8}[lemma65]{Lemma}
\begin{lemma8}\label{lemma8}
Assume $\phi\in C^\infty_0(\mathbb{R}^d)$, $\phi(0)=1$, and
$\phi(x)=0$ if $|x|\geq 1$. Then for any $\lambda>0$ and $g\in X$,
\begin{align*}
&\Im \langle g,R_0(\lambda\pm i0) g\rangle =c_1 \int_{\sqrt{\lambda} S^{d-1}} |\hat{g}(\xi)|^2 \, \sigma(d\xi)\\
&\lim_{R\to\infty}\int_{\mathbb{R}^d} |R_0(\lambda\pm i0)
g|^2(x)\, \phi\big(\frac{x}{R}\big)\, \frac{dx}{R}=c_2
\int_{\sqrt{\lambda}S^{d-1}} |\hat{g}(\xi)|^2\, \sigma(d\xi)
\end{align*}
where $c_1=c_1(\lambda,\pm)\neq 0$ and $
c_2=c_2(\lambda,\pm,\phi)\neq 0$.
\end{lemma8}
\begin{proof}[Proof of Lemma~\ref{lemma8}]
If $g\in\mathcal{S}(\mathbb{R}^d)$, then these properties are
standard, see for example \cite[Chapter XIV]{Ho}. Clearly,
$\mathcal{S}(\mathbb{R}^d)$ is dense in~$X$. Moreover, in view of
Theorem~\ref{main2} the left-hand sides of these limits are
continuous with respect to the norm of $X$. Finally, by the
trace lemma and the Stein-Tomas theorem, respectively, the
right-hand sides are also continuous with respect to the
$X$-norm, which proves the identities.
\end{proof}

Assume from now on that $L$ is the admissible perturbation in Theorem \ref{main1}. We define a set $\widetilde{\mathcal E}\subset\mathbb
R\setminus\{0\}$ so that off this set~\eqref{eq:fred} is
invertible. We will show later that $\widetilde{\mathcal E}$ is
exactly the set of nonzero eigenvalues, which we denoted by
$\mathcal{E}$ in Theorem~\ref{main1}. Let
\begin{equation}\label{hu1}
\widetilde{\mathcal{E}}^{\pm}:=\{\lambda\in\mathbb{R}\setminus\{0\}:
\text{there is }f\in X^\ast\setminus\{0\}\text{ with
}(\mathrm{Id}_{X^\ast}+R_0(\lambda\pm i0)L)f=0\}.
\end{equation}
Notice that
$R_0(\lambda-i0)\overline{g}=\overline{R_0(\lambda+i0)g}$ for any
$g\in X$. Thus
$\widetilde{\mathcal{E}}^{+}=\widetilde{\mathcal{E}}^{-}:=\widetilde{\mathcal{E}}$.
For any $\lambda\in\mathbb{R}\setminus \{0\}$ we define the
eigenspaces
\begin{equation*}
\mathcal{F}_\lambda^{\pm}=\{f\in
X^\ast:(\mathrm{Id}_{X^\ast}+R_0(\lambda\pm i0)L)f=0\}.
\end{equation*}

\newtheorem{lemma9.1}[lemma65]{Lemma}
\begin{lemma9.1}\label{lemma9.1}
Assume that $\lambda\in\widetilde{\mathcal{E}}$ and
$f\in\mathcal{F}_\lambda^{\pm}$. Then $(1+|x|^2)^N f\in X^\ast$ for any $N\geq 0$, and
\begin{equation}
\label{eq:wei}
 ||(1+|x|^2)^Nf||_{X^\ast} \le C_{N,L,\lambda} \|f\|_{X^\ast}.
\end{equation}
\end{lemma9.1}

\begin{proof}[Proof of Lemma~\ref{lemma9.1}] We show first that
\begin{equation}\label{tb600}
\langle Lf,g\rangle=\overline{\langle Lg,f\rangle}
\end{equation}
for any $f,g\in X^\ast$ (this is assumed in \eqref{tb61} for $f,g\in\mathcal{S}(\mathbb{R}^d)$). Using the functions $\chi$ and $\varphi$ defined at the end of section \ref{weight} we define the sequences
\begin{equation*}
f_n(x)=\chi_n(x)(f\ast\varphi_{1/n})(x)\text{ and }g_{n}(x)=\chi_n(x)(g\ast\varphi_{1/n})(x).
\end{equation*}
Clearly, $f_n,g_n\in\mathcal{S}(\mathbb{R}^d)$, $||f_n||_{X^\ast}\leq C||f||_{X^\ast}$, and $||g_n||_{X^\ast}\leq C||g||_{X^\ast}$. In view of \eqref{tb61}, it suffices to prove that
\begin{equation}\label{tb660}
\lim_{n\to\infty}\langle Lf_n,g_n\rangle=\langle Lf,g\rangle.
\end{equation}
We remark that the sequences $f_n$ and $g_n$ may not converge to $f$ and $g$ in $X^\ast$ (in fact $\mathcal{S}(\mathbb{R}^d)$ is not dense in $X^\ast$). However, using \eqref{tb5} and the fact that $\mathcal{S}(\mathbb{R}^d)$ is dense in $X$, for the limit above it suffices to prove that
\begin{equation}\label{tb601}
\lim_{n\to\infty}||L(f_n-f)||_{X}=0,
\end{equation}
and
\begin{equation}\label{tb602}
\lim_{n\to\infty}\langle(g_n-g),\phi\rangle=0\text{ for any }\phi\in\mathcal{S}(\mathbb{R}^d).
\end{equation}
The limit \eqref{tb602} is clear, even for any $g\in L^{p'_d}$. For the limit \eqref{tb601}, given $\varepsilon>0$, we use \eqref{tb6} with $N=0$. The result is
\begin{equation*}
||L(f_n-f)||_{X}\leq \varepsilon||f_n-f||_{X^\ast}+A_\varepsilon||f_n-f||_{L^2(\{|x|\leq R_\varepsilon\})}.
\end{equation*}
Since $\lim_{n\to\infty}||f_n-f||_{L^2(\{|x|\leq R\})}=0$ for any $R\geq 1$, the limit \eqref{tb601} follows.

Assume that
\begin{equation}
\label{eq:fdef}
 f + R_0(\lambda\pm i0)L f =0
\end{equation}
for some $f\in X^*$, $\lambda\in\mathbb{R}\setminus\{0\}$, and
some choice of $+$ or $-$. If $\lambda>0$, then we use
Lemma~\ref{lemma8} and the fact that $\langle
Lf,f\rangle\in\mathbb{R}$ to conclude that
\[ 0=\Im \langle Lf,R_0(\lambda\pm i0)Lf\rangle = c \int_{\sqrt{\lambda} S^{d-1}} |\widehat{Lf}|^2 \, d\sigma\]
with some constant $c\ne0$. Applying Lemma~\ref{lemma8} again
implies that
\[ \lim_{R\to\infty} R^{-1} \int_{\{|x|\le R\}} \Big| [R_0(\lambda\pm i0)Lf] (x)\Big|^2 \, dx =0.\]
In view of~\eqref{eq:fdef} this is the same as
\begin{equation}\label{new1}
\lim_{R\to\infty} R^{-1} \int_{\{|x|\le R\}} |f (x)|^2 \, dx =0.
\end{equation}
Let $\delta>0$ be such that $\delta\le |\lambda|\le \delta^{-1}$.
By Lemma~\ref{lemma65} and Theorem~\ref{main3}
\[ ||\mu_{N,\gamma}f||_{X^\ast}\leq C_{N,\delta}||\mu_{N,\gamma}Lf||_{X}.\]
We use \eqref{tb6} with $\varepsilon=(2C_{N,\delta})^{-1}$ and the
fact that $||\mu_{N,\gamma}f||_{X^\ast}<\infty$. By absorbing the
term $(1/2)||\mu_{N,\gamma}f||_{X^\ast}$,
\begin{equation*}
||\mu_{N,\gamma}f||_{X^\ast}\leq C_{N,L,\delta}||f||_{B^\ast}.
\end{equation*}
The inequality \eqref{eq:wei} follows by letting $\gamma\to 0$.

If $\lambda<0$, then since $Lf\in X\hookrightarrow W^{-1,2}$, we
have $R_0(\lambda)L f\in W^{1,2}$, thus \eqref{new1} holds, using
\eqref{eq:fdef}. The same argument as above proves~\eqref{eq:wei}.
\end{proof}

\newtheorem{lemma9}[lemma65]{Lemma}
\begin{lemma9}\label{lemma9}
(a) For any $\lambda\in\mathbb{R}\setminus\{0\}$
\begin{equation*}
\mathcal{F}_\lambda^{\pm}\subseteq\{u\in\mathrm{Domain}(H):Hu=\lambda
u\}.
\end{equation*}
In particular, $\widetilde{\mathcal{E}}\subseteq\mathcal{E}$.

(b) The set $\widetilde{\mathcal{E}}$ is discrete in
$\mathbb{R}\setminus\{0\}$, i.e., $I\cap\widetilde{\mathcal{E}}$
is finite for any compact set
$I\subseteq\mathbb{R}\setminus\{0\}$.

(c) For any $\lambda\in\mathbb{R}\setminus\{0\}$, the vector
spaces $\mathcal{F}_\lambda^{\pm}$ are finite-dimensional.
\end{lemma9}

\begin{proof}[Proof of Lemma~\ref{lemma9}]
For part (a), by Lemma~\ref{lemma65} (c)
\begin{equation}\label{hu6}
(-\Delta-\lambda)f+Lf=0
\end{equation}
for any $f\in\mathcal{F}_\lambda^{\pm}$. By Lemma~\ref{lemma9.1},
$f\in W^{1,2}$, thus $f\in\mathrm{Domain(H)}$.

We prove now part (b). Assume, for contradiction, that the set
$\widetilde{\mathcal{E}}\cap\{\lambda:\delta\leq|\lambda|\leq\delta^{-1}\}$
is infinite for some $\delta>0$, thus
$\widetilde{\mathcal{E}}\cap\{\lambda:\delta\leq|\lambda|\leq\delta^{-1}\}=\{\lambda_1,\lambda_2,\ldots\}$,
$\lambda_m\neq\lambda_n$ if $m\neq n$. For any $n$ fix $f_n\in
\mathcal{F}_{\lambda_n}^+\setminus\{0\}$. By \eqref{hu6}, $f_m\neq
f_n$ if $m\neq n$. By \eqref{eq:wei}, $f_n\in W^{1,2}$. We
normalize the functions $f_n$ in such a way that
$||f_n||_{W^{1,2}}=1$. Then, by \eqref{eq:wei},
\begin{equation}\label{hu8}
||(1+|x|^2)f_n||_{W^{1,2}}\leq C_{L,\delta}||f_n||_{X^\ast}\leq
C_{L,\delta}
\end{equation}
for any integer $n\geq 1$. Also, by \eqref{hu6},
$(-\Delta+1)f_n=(\lambda_n+1-L)f_n$, thus
\begin{equation}\label{hu10}
\begin{split}
f_n&=R_0(-1)[(\lambda_n+1-L)f_n]\\
&=(\lambda_n+1)R_0(-1)(1+|x|^2)^{-1}[(1+|x|^2)f_n]-R_0(-1)Lf_n.
\end{split}
\end{equation}
Using Lemma~\ref{lemma7}, it is easy to see that the operators
\begin{equation}\label{hu7}
R_0(-1)(1+|x|^2)^{-1},R_0(-1)L: W^{1,2}\to W^{1,2}\text{ as
compact operators}.
\end{equation}
By \eqref{hu8} and \eqref{hu10}, we pass to a subsequence and
assume that $\lim_{n\to\infty}f_n=f_\infty$ in $W^{1,2}$. By the
normalization of the functions $f_n$, $||f_\infty||_{W^{1,2}}=1$.
On the other hand, by part (a), the functions $f_n$ are
eigenfunctions of the self-adjoint operator $H$ (see section
\ref{proof}) with different eigenvalues, thus $\langle
f_m,f_n\rangle=0$ if $m\neq n$. Therefore $f_\infty=0$, which
yields a contradiction.

For part (c), assume for contradiction that
$\mathrm{dim}(\mathcal{F}_\lambda^{\pm})=\infty$ for some
$\lambda\in\mathbb{R}\setminus\{0\}$. Then we could find an
infinite sequence of functions $\{f_n\}_{n=1}^\infty\subseteq
\mathcal{F}_\lambda^{\pm}$, such that $||f_n||_{W^{1,2}}=1$ and
$\langle f_m,f_n\rangle=0$ if $m\neq n$. The same argument as
above gives a contradiction.
\end{proof}

Finally, we need to prove that the operators
$\mathrm{Id}_{X^\ast}+R_0(\lambda\pm i\epsilon)L$ are uniformly
invertible on $X^\ast$, provided that $\lambda$ is separated
from $\widetilde{\mathcal{E}}$.

\newtheorem{lemma10}[lemma65]{Lemma}
\begin{lemma10}\label{lemma10}
(a) For any $\lambda\in(\mathbb{R}\setminus\{0\})\setminus
\widetilde{\mathcal{E}}$ and $\epsilon\in[0,\infty)$, the operator
$\mathrm{Id}_{X^\ast}+R_0(\lambda\pm i\epsilon) L$ is invertible
on $X^\ast$.

(b) For any compact set $I\subseteq
(\mathbb{R}\setminus\{0\})\setminus \widetilde{\mathcal{E}}$,
\begin{equation}\label{hu20}
\sup_{\lambda\in I} \sup_{1\ge\epsilon\ge0} \Big\|
(\mathrm{Id}_{X^\ast}+R_0(\lambda\pm i\epsilon)L)^{-1}
\Big\|_{X^\ast\to X^\ast} < \infty.
\end{equation}
\end{lemma10}
\begin{proof}[Proof of Lemma~\ref{lemma10}] For part (a) we use
Lemma~\ref{lemma7}. Since $R_0(\lambda\pm i\epsilon)L$ is compact
on $X^\ast$, the only alternative to invertibility is the existence of a
nontrivial kernel. By the definition of the set $\widetilde{E}$,
such a nontrivial kernel could only exist if $\epsilon>0$. If $f\in
X^\ast$ has the property that
\begin{equation}\label{hu30}
f+R_0(\lambda\pm i\epsilon)Lf=0,
\end{equation}
then
\begin{equation*}
\langle Lf,f\rangle+\langle Lf,R_0(\lambda+i\epsilon)Lf\rangle=0.
\end{equation*}
Since $L$ is symmetric, by taking the imaginary part we have
\begin{equation*}
\begin{split}
0=\Im\langle
Lf,R_0(\lambda+i\epsilon)Lf\rangle=\epsilon\int_{\mathbb{R}^d}|\widehat{Lf}|^2
[(|\xi|^2-\lambda)^2+\epsilon^2]^{-1}\,d\xi.
\end{split}
\end{equation*}
Since $\epsilon\neq 0$, it follows that $Lf\equiv 0$, thus
$f\equiv 0$ by \eqref{hu30}.

For part (b), we show that
\begin{equation}\label{hu201}
\sup_{\lambda\in I} \sup_{1\ge\epsilon\ge0} \Big\|
(\mathrm{Id}_{X^\ast}+R_0(\lambda+ i\epsilon)L)^{-1}
\Big\|_{X^\ast\to X^\ast} < \infty.
\end{equation}
The proof for the operators $(\mathrm{Id}_{X^\ast}+R_0(\lambda-i\epsilon)L)^{-1}$ is identical. Assume for contradiction that the supremum on the
left-hand side of \eqref{hu201} is $\infty$. Then there exist
$f_n\in X^\ast$, $\|f_n\|_{X^*} =1$, such that
\[ \|(\mathrm{Id}_{X^\ast}+R_0(\lambda_n+i\epsilon_n)L) f_n \|_{X^*} \to 0 \]
as $n\to\infty$. Here $\lambda_n\in I$ and $\epsilon_n\in[0,1]$. We start from the identity
\begin{equation}\label{tb80}
f_n=-R_0(\lambda_n+i\epsilon_n)Lf_n+r_n,
\end{equation}
where $||r_n||_{X^\ast}\to 0$ as $n\to\infty$. By passing to a
subsequence, we may assume that $\lambda_n+i\epsilon_n\to \lambda_\infty+i\epsilon_\infty\in I\times[0,1]$ and, since $L:X^\ast\to X$ is compact (Lemma
\ref{lemma7}), $Lf_n\to h$ in X. Let $f_\infty=-R_0(\lambda_\infty+i\epsilon_\infty)h\in X^\ast$. Using \eqref{tb80}, for any $\phi\in\mathcal{S}(\mathbb{R}^d)$ we have
\begin{equation}\label{pp60}
\begin{split}
&\lim_{n\to\infty}\langle f_n,\phi\rangle=\lim_{n\to\infty}\langle -R_0(\lambda_n+i\epsilon_n)Lf_n,\phi\rangle\\
&=\lim_{n\to\infty}\langle -R_0(\lambda_n+i\epsilon_n)h,\phi\rangle+\lim_{n\to\infty}\langle R_0(\lambda_n+i\epsilon_n)(h-Lf_n),\phi\rangle=\langle f_\infty,\phi\rangle.
\end{split}
\end{equation}
In the last identity we used Lemma \ref{lemma65} and Theorem \ref{main2}. In addition, for any $\varepsilon>0$ we use \eqref{tb6} with $N=0$:
\begin{equation*}
\begin{split}
||Lf_n-Lf_\infty||_{X}&\leq \varepsilon||f_n-f_\infty||_{X^\ast}+A_\varepsilon||(f_n-f_\infty)\mathbf{1}_{\{|x|\leq R_\varepsilon\}}||_{L^2}\\
&\leq \varepsilon||f_n-f_\infty||_{X^\ast}+C_\varepsilon[||r_n||_{X^\ast}+||R_0(\lambda_n+i\epsilon_n)(Lf_n-h)||_{X^\ast}]\\
&+A_\varepsilon||\mathbf{1}_{\{|x|\leq R_\varepsilon\}}[R_0(\lambda_n+i\epsilon_n)h-R_0(\lambda_\infty+i\epsilon_\infty)h]||_{L^2}.
\end{split}
\end{equation*}
By Lemma \ref{lemma65} and Theorem \ref{main2}
\begin{equation}\label{pp61}
\lim_{n\to\infty}Lf_n=Lf_\infty\text{ in }X.
\end{equation}
It follows from \eqref{pp60}, \eqref{pp61} and Lemma \ref{lemma65} that for any $\phi\in\mathcal{S}(\mathbb{R}^d)$
\begin{equation*}
0=\lim_{n\to\infty}\langle f_n+R_0(\lambda_n+i\epsilon_n)Lf_n,\phi\rangle=\langle f_\infty+R_0(\lambda_\infty+i\epsilon_\infty)Lf_\infty,\phi\rangle.
\end{equation*}
Thus $f_\infty+R_0(\lambda_\infty+i\epsilon_\infty)Lf_\infty=0$, which, in view of part (a), shows that $f_\infty=0$. By \eqref{pp61}, $\lim_{n\to\infty}Lf_n=0$ in $X$. This gives a contradiction, in view of the identity \eqref{tb80} and the fact that $||f_n||_{X^\ast}=1$.
\end{proof}

\section{Proof of Theorem \ref{main1}}\label{proof}

We can now finish the proof of Theorem \ref{main1}.

{\it{Proof of part (a):}} We show
first that if $u\in\mathrm{Domain}(H)$ then
\begin{equation}\label{gu2}
||u||_{W^{1,2}}^2\leq C_L||u||_{L^2}^2+C_0\langle Hu,u\rangle,
\end{equation}
for some constants $C_0\geq 0$ and $C_L$. To see this, we start from the identity
\begin{equation*}
\langle Hu,\phi\rangle=\langle\nabla u,\nabla\phi\rangle+\langle Lu,\phi\rangle,
\end{equation*}
valid, by definition, for any $u\in \mathrm{Domain}(H)$ and
$\phi\in\mathcal{S}(\mathbb{R}^d)$. Since
$u\in\mathrm{Domain}(H)$, and $\mathcal{S}(\mathbb{R}^d)$ is dense
in $W^{1,2}(\mathbb{R}^d)$, we conclude that
\begin{equation}\label{gu3}
\langle Hu,v\rangle=\langle\nabla u,\nabla v\rangle+\langle Lu,v\rangle,
\end{equation}
for any $u,v\in \mathrm{Domain}(H)$. Here we used that $L:W^{1,2}\to W^{-1,2}$.
In particular, $\langle Hu,u\rangle\in\mathbb{R}$ for any $u\in \mathrm{Domain}(H)$.
For $\varepsilon>0$ small enough, we use \eqref{tb6} with $N=0$. The result is
\begin{equation*}
\begin{split}
\langle Hu,u\rangle&\geq ||\,|\nabla u|\,||_{L^2}^2-|\langle Lu,u\rangle|\\
&\geq c_0||u||_{W^{1,2}}^2-C||u||_{L^2}^2-||u||_{X^\ast}(\varepsilon||u||_{X^\ast}+A_\varepsilon||u||_{L^2})\\
&\geq c_0||u||_{W^{1,2}}^2-C||u||_{L^2}^2-C\varepsilon||u||_{W^{1,2}}^2-C_\varepsilon||u||_{L^2}^2\\
&\geq (c_0/2)||u||_{W^{1,2}}^2-C_{L}||u||_{L^2}^2,
\end{split}
\end{equation*}
by choosing $\varepsilon$ small enough. This proves \eqref{gu2}.

An elementary limiting argument, using \eqref{gu2} and the fact
that $W^{1,2}$ is a Banach space, shows that the
operator
\begin{equation*}
H:\mathrm{Domain}(H)\to L^2
\end{equation*}
is closed. Clearly, it is also bounded from below. Finally, using
\eqref{gu3} and the fact that $L$ is symmetric, the operator $H$
is symmetric. It remains to prove that
$\mathrm{Domain}(H)$ is dense in $L^2$, and, using the criterion
for self-adjointness \cite[Theorem VIII.3]{RS2}, that
\begin{equation}\label{gu4}
\mathrm{Range}(H\pm i)=L^2.
\end{equation}

\newtheorem{lemma99}{Lemma}[section]
\begin{lemma99}\label{lemma99}
If $\lambda\in\mathbb{R}$, and
$\epsilon\in\mathbb{R}\setminus\{0\}$, then
\begin{equation*}
\widetilde{R}_L(\lambda+i\epsilon):=(\mathrm{Id}_{W^{1,2}}+R_0(\lambda+i\epsilon)L)^{-1}R_0(\lambda+i\epsilon)
\end{equation*}
defines a bounded operator
$\widetilde{R}_L(\lambda+i\epsilon):L^2\to\mathrm{Domain}(H)$. In addition,
\begin{equation}\label{gu6}
[H-(\lambda+i\epsilon)]\widetilde{R}_L(\lambda+i\epsilon)=\mathrm{Id}_{L^2}.
\end{equation}
\end{lemma99}

Assuming Lemma~\ref{lemma99}, the density of $\mathrm{Domain}(H)$
in $L^2$ and \eqref{gu4} follow easily (for the density of
$\mathrm{Domain}(H)$, notice that $\mathrm{Domain}(H)$ contains
the set $(\mathrm{Id}_{W^{1,2}}+R_0(i)L)^{-1}W^{2,2}$, which is
dense in $W^{1,2}$, thus dense in $L^2$). In addition, the
identity \eqref{gu3} shows that the map
$H-(\lambda+i\epsilon):\mathrm{Domain}(H)\to L^2$, $\epsilon\neq
0$, is injective. Thus, the spectrum of the operator $H$ is a
subset of $\mathbb{R}$, and we have the resolvent identity
\begin{equation}\label{gu8}
R_L(\lambda+i\epsilon)=(\mathrm{Id}_{W^{1,2}}+R_0(\lambda+i\epsilon)L)^{-1}R_0(\lambda+i\epsilon),
\end{equation}
for $\lambda\in\mathbb{R}$ and $\epsilon\in\mathbb{R}\setminus\{0\}$.

\begin{proof}[Proof of Lemma~\ref{lemma99}] We show first that the operator
$(\mathrm{Id}_{W^{1,2}}+R_0(\lambda+i\epsilon)L)$ is well-defined and invertible on $W^{1,2}$.
Using \eqref{tb2}, Lemma~\ref{lemma7}, and the fact that $\epsilon\neq 0$,
the operator $R_0(\lambda+i\epsilon)L:W^{1,2}\to W^{1,2}$ is bounded and compact. By Fredholm's
alternative, it suffices to prove that the kernel of this operator
is trivial. Assume $f\in W^{1,2}$ has the property that
\begin{equation*}
f+R_0(\lambda+i\epsilon)Lf=0.
\end{equation*}
The same argument as in the proof of Lemma~\ref{lemma10}(a) shows that $f\equiv 0$,
which completes the proof of invertibility.

Therefore, the map $\widetilde{R}_L(\lambda+i\epsilon):L^2\to
W^{1,2}$ is a bounded operator. It remains to verify the identity
\eqref{gu6}. Assume $f\in L^2$ and let
$g=R_0(\lambda+i\epsilon)f\in W^{2,2}$ and
$h=(\mathrm{Id}_{W^{1,2}}+R_0(\lambda+i\epsilon)L)^{-1}g\in
W^{1,2}$. Then
\begin{equation*}
h=g-R_0(\lambda+i\epsilon)Lh.
\end{equation*}
Thus, in $\mathcal{S}'(\mathbb{R}^d)$ we have
\begin{equation*}
\begin{split}
[-\Delta+L-(\lambda+&i\epsilon)]h=Lh+[-\Delta-(\lambda+i\epsilon)]g\\
&-[-\Delta-(\lambda+i\epsilon)]R_0(\lambda+i\epsilon)Lh=f,
\end{split}
\end{equation*}
as desired.
\end{proof}

\medskip

{\it{Proof of part (c):}} Assume $u\in\mathrm{Domain}(H)$ and
$Hu=\lambda u$, $\lambda\in\mathbb{R}\setminus\{0\}$. Since $u\in
W^{1,2}$, we have $(\Delta+\lambda)u=Lu$, $u\in X^\ast$, and
Theorem~\ref{main3} applies. Thus
\[ ||\mu_{N,\gamma}u||_{X^\ast}\leq C_{N,\lambda}||\mu_{N,\gamma}Lu||_{X}.\]
As in the proof of Lemma~\ref{lemma9.1} we use \eqref{tb6} with
$\varepsilon=(2C_{N,\lambda})^{-1}$ and the fact that
$||\mu_{N,\gamma}u||_{X^\ast}<\infty$. By absorbing the term
$(1/2)||\mu_{N,\gamma}u||_{X^\ast}$,
\begin{equation*}
||\mu_{N,\gamma}u||_{X^\ast}\leq C_{N,L,\lambda}||u||_{B^\ast}.
\end{equation*}
Part (iii) follows by letting $\gamma\to 0$.
\medskip

{\it{Proof of part (b):}} For any $\lambda\in\mathcal{E}$ let
\begin{equation*}
\mathcal{H}_\lambda=\{u\in\mathrm{Domain}(H):Hu=\lambda u\}.
\end{equation*}
By Lemma~\ref{lemma9},
$\widetilde{\mathcal{E}}\subseteq\mathcal{E}$ and
$\mathcal{F}_\lambda^{+}\cup\mathcal{F}_\lambda^{-}\subseteq\mathcal{H}_\lambda$.
It suffices to show that
$\mathcal{E}\subseteq\widetilde{\mathcal{E}}$ and
$\mathcal{H}_\lambda\subseteq\mathcal{F}_\lambda^{+}\cap\mathcal{F}_\lambda^{-}$.
Since $\mathrm{Domain}(H)\subseteq X^\ast$, it suffices to show
that if $u\in\mathcal{H}_\lambda$ then
\begin{equation}\label{ga1}
u+R_0(\lambda\pm i0)Lu=0.
\end{equation}
Since $(-\Delta-\lambda)u+Lu=0$ we have
\begin{equation*}
R_0(\lambda\pm i0)[(-\Delta-\lambda)u]+R_0(\lambda\pm i0)Lu=0.
\end{equation*}
For \eqref{ga1}, it suffices to prove that
\begin{equation*}
R_0(\lambda\pm i0)[(-\Delta-\lambda)u]=u.
\end{equation*}
This is clear if $\lambda<0$, for any $u\in\mathcal{S}'(\mathbb{R}^d)$. Assume $\lambda>0$ and $R_0(\lambda\pm
i0)[(-\Delta-\lambda)u]=u'\in X^\ast$. By Lemma~\ref{lemma65}(c),
since $(-\Delta-\lambda)u=-Lu\subseteq X$,
$(-\Delta-\lambda)u'=(-\Delta-\lambda)u$, thus
\begin{equation}\label{ga2}
(-\Delta-\lambda)(u-u')=0.
\end{equation}
Since $(-\Delta-\lambda)u\in X$, and by definition of $u'$,
 Lemma~\ref{lemma8} gives
\begin{equation}
\label{eq:up}
\lim_{R\to\infty}\int_{\mathbb{R}^d} |u'|^2(x)\, \phi(x/R)\,dx/R=c_2(\lambda)
\int_{\sqrt{\lambda}S^{d-1}} |\widehat{(-\Delta-\lambda)u}|^2\, \sigma(d\xi).
\end{equation}
Since $u$ is rapidly decreasing in $L^2$ (using part (c)), it follows that $u\in L^1(\mathbb{R}^d)$,
and thus $\hat{u}\in C(\mathbb{R}^d)$. Hence,
\[  [(-\Delta-\lambda)u]^{\wedge}(\xi) = (\xi^2-\lambda)\hat{u}(\xi)\]
both in the sense of distributions and as continuous functions. But the right-hand side vanishes
on $\sqrt{\lambda} S^{d-1}$, and so the limit in~\eqref{eq:up} vanishes.
Thus
\begin{equation*}
\lim_{R\to\infty}R^{-1}\int_{|x|\leq R}|u'|^2\,dx=0.
\end{equation*}
Using again the fact that $u$ is rapidly decreasing in $L^2$, we can apply Theorem~\ref{main3} with $N=0$ to the function $u-u'$. The identity \eqref{ga2} gives $u\equiv u'$, which completes the proof of \eqref{ga1}.
\medskip

{\it{Proof of part (d):}} We use the resolvent identity
\begin{equation*}
R_L(\lambda+i\epsilon)=(\mathrm{Id}_{W^{1,2}}+R_0(\lambda+i\epsilon)L)^{-1}R_0(\lambda+i\epsilon)
\end{equation*}
(see \eqref{gu8}). Recall that
$\mathcal{E}=\widetilde{\mathcal{E}}$. The main bound \eqref{bi3}
then follows from Theorem~\ref{main2} and Lemma~\ref{lemma10}(b).
\medskip

{\it{Proof of part (e):}} The statement
$\sigma_{\mathrm{sc}}(H)=\emptyset$ is an immediate consequence of
part (d), see~\cite[Theorem XIII.20]{RS4}. To prove that
$\sigma_{\mathrm{ac}}(H)\subseteq [0,\infty)$, assume
$\lambda\in(-\infty,0)\setminus\mathcal{E}$. We have to prove that
$\lambda$ is in the resolvent set of $H$. We use
Lemma~\ref{lemma99}. Since $\lambda\notin\sigma_{\mathrm{pp}}(H)$,
the equation $(\mathrm{Id}_{W^{1,2}}+R_0(\lambda)L)f=0$ has no
solutions in $W^{1,2}$. Since the operator $R_0(\lambda)L$ is
compact on $W^{1,2}$ (see the proof of Lemma~\ref{lemma99}),
Fredholm's alternative shows that the operator
$(\mathrm{Id}_{W^{1,2}}+R_0(\lambda)L)$ is invertible on
$W^{1,2}$. It follows, as in Lemma~\ref{lemma99}, that $\lambda$
is in the resolvent set of $H$.

The reverse inclusion
$\sigma_{\mathrm{ac}}(H)\supseteq [0,\infty)$ follows from the
existence of the wave operators which we establish in the next
paragraph.
\medskip

{\it{Proof of part (f):}}  This will be done by means of a local version of
Kato's smoothing theory. This is the only place in the proof where condition~(3)
of our definition of admissible perturbations is required.
According to~\cite[Theorem XIII.31]{RS4} and its corollary\footnote{
Strictly speaking, \cite[Theorem XIII.31]{RS4} and its corollary are only stated with~$J=1$.
But the same proof also applies to the case $J>1$ needed here. Indeed, the only change is
to the first inequality on page~166 of~\cite{RS4} which  needs to be replaced with
\[ \le \sum_{j=0}^J \|A_j(H-z)^{-1}\|\, \|B_j(H_0-z)^{-1} e^{-iH_0t} E_I^{(0)} \phi\|\, \|\psi\|.\]
},
we need to prove the following: Write $H_0=-\Delta$, $H-H_0=L=\sum_{j=1}^J A_j^*B_j$
as in condition~(3). Then we need to show that each $B_j$ is $H_0$-bounded
and that each $A_j$ is $H$-bounded. Furthermore, we need to show that for any compact interval $I$
so that
\begin{equation}
 \label{eq:I}
 I\subset \mathbb{R}\setminus (\mathcal{E}\cup\{0\})
\end{equation}
we have the property that $A_j E(I)$ is $H$-smooth and $B_j E_0(I)$ is $H_0$-smooth
in the sense of Kato, see~\cite{kato}.
Here $E_0$ and $E$ denote the spectral projections associated with $H_0$ and $H$,
respectively.

We start with the boundedness properties, and then discuss the smoothness.
Thus we need to prove that there exist constants $a,b$ so that for each $1\le j\le J$
\begin{align}
\mathrm{Domain}(A_j) & \supseteq \mathrm{Domain}(-\Delta) \nonumber\\
\|A_j f\|_{L^2} &\le a\|\Delta f\|_{L^2} + b\|f\|_{L^2} \qquad \forall \; f\in \mathrm{Domain}(\Delta) \label{eq:Abd}\\
\mathrm{Domain}(B_j) & \supseteq \mathrm{Domain}(H) \nonumber \\
\|B_j f\|_{L^2} &\le a\|H f\|_{L^2} + b\|f\|_{L^2} \qquad \forall
\; f\in \mathrm{Domain}(H). \label{eq:Bbd}
\end{align}
By assumption, $\mathrm{Domain}(A_j),\mathrm{Domain}(B_j)\supseteq
W^{1,2}$, so that the required set inclusions are clear.
Furthermore, condition~(3) guarantees that
\begin{equation*}
\|A_j f\|_{L^2} + \|B_j f\|_{L^2} \le C\|f\|_{X^\ast}\leq
C\|f\|_{W^{1,2}}.
\end{equation*}
In conjunction with~\eqref{gu2}, this implies~\eqref{eq:Abd} and~\eqref{eq:Bbd}.

Next, we discuss the smoothness of $A_j$ and $B_j$. In view
of~\eqref{eq:I} the limiting absorption principle~\eqref{bi3}
holds for~$I$, and similarly for the free resolvent $R_0$. It is
shown in \cite[Theorem~XIII.30]{RS4} that it suffices to prove
that
\[ \sup_{\lambda\in I,0<\epsilon<1} |\eps|\,\big\|\, A_j R_L(\lambda+i\epsilon) \,
\big\|_{L^2\to L^2}^2 < \infty,\,\sup_{\lambda\in I,0<\epsilon<1}
|\eps|\,\big\|\, B_j R_0(\lambda+i\epsilon)\, \big\|_{L^2\to
L^2}^2 < \infty
\]
for the required smoothness properties of $A_j$ and $B_j$ to hold.
However, these are easy consequences of~\eqref{bi3}
and~\eqref{bu1}, respectively, since we are requiring that
$A_j,B_j:X^*\to L^2$ as bounded operators. Indeed, we only need to
verify that \begin{equation} \label{eq:2X*}
\sup_{\substack{0<\epsilon<1\\
\lambda\in I}} |\eps|\,\big\| R_L(\lambda+i\epsilon) \,
\big\|_{L^2\to X^*}^2 < \infty, \quad \sup_{\substack{0<\epsilon<1\\
\lambda\in I}} |\eps|\,\big\|R_0(\lambda+i\epsilon)\,
\big\|_{L^2\to X^*}^2 <\infty. \end{equation} To see this, fix
$\eps\ne0$ and apply the resolvent indentity with $f\in X$:
\begin{align*}
\|R_L(\lambda+ i\eps)f\|_{L^2}^2 &= \langle
R_L(\lambda+ i\eps)^*R_L(\lambda+i\eps)f,f\rangle \\
&= \frac{-1}{2i\eps} \Big\langle \big(
R_L(\lambda+i\eps)^*-R_L(\lambda+i\eps)\big)f,f \Big\rangle \\
&\le \frac{1}{2|\eps|} [\|R_L(\lambda+ i\eps)f\|_{X^*}+
\|R_L(\lambda- i\eps)f\|_{X^*}]\|f\|_X \\
& \le \frac{C}{|\eps|} \|f\|_X^2,
\end{align*}
by \eqref{bi3}, and similarly for $R_0(\lambda+i\eps)$.  Now
suppose $g\in L^2$ and $f\in X$. Then this estimate implies that
\[ \big|\langle R_L(\lambda+i\eps)g, f\rangle\big| \le
C|\eps|^{-1/2}\, \|f\|_X \|g\|_{L^2}.
\]
Thus, $R_L(\lambda+i\eps)g$ is an element of $X^*$ with norm
\[ |\eps|\,\|R_L(\lambda+i\eps)g\|_{X^*}^2 \le C\|g\|^2_{2},\]
and similarly for $R_0(\lambda+i\eps)g$. Hence, we are done, i.e.,
the wave operators
\[ \Omega^{\pm}(H,H_0):=\slim_{t\to\mp \infty} e^{itH} e^{-itH_0},
\quad \Omega^{\pm}(H_0,H):=\slim_{t\to\mp \infty} e^{itH_0}
e^{-itH}E_{a.c} \] exist and are complete, see the aforementioned
corollary in~\cite{RS4}.

We now return to the issue of showing that $\sigma(H)\cap I\ne\emptyset$
for any nonempty interval $I\subset [0,\infty)$. Indeed, fix any such compact interval
which also satisfies~\eqref{eq:I} and let $W_{\pm}, \tilde{W}_{\pm}$ be the local wave operators
defined as the strong limits
\begin{equation}
\label{eq:wop} W_{\pm}:=\slim_{t\to\mp\infty} e^{iHt} e^{-itH_0}
E_0(I), \quad \tilde{W}_{\pm}:=\slim_{t\to\mp\infty} e^{iH_0t}
e^{-itH} E(I).
\end{equation}
These strong limits exist because of \cite[Theorem~XIII.31]{RS4}.
Moreover, the relations
\[ W_{\pm}^*=\tilde{W}_{\pm},\quad \tilde{W}_{\pm}W_{\pm}=E_0(I),\quad
W_{\pm}\tilde{W}_{\pm}=E(I)
\]
hold. Since $E_0(I)\ne0$ by choice of $I$, it
follows that $W_{\pm}$ is an isometry on the range of $E_0(I)$ and $\|W_{\pm}\|=1$.
Thus also $\|W_{\pm}^*\|=1$. Choose any $f\in L^2$ with $W^*_{\pm}f\ne0$ and observe
that
\[ \|E(I)f\|_{L^2}^2 = \langle W_{\pm}W_{\pm}^* f,f\rangle = \|W_{\pm}^* f\|_{L^2}^2 \ne 0.\]
Hence $E(I)\ne0$, which shows that $\sigma(H)\cap I\ne\emptyset$, as claimed.

{\it{Proof of Corollary \ref{mai2}:}}
Note that any  $F\in L^{d+1}(\mathbb{R}^d)$ satisfies, by Sobolev imbedding,
\[ \|FS_{1/(d+1)} f\|_{L^2} \le C\|F\|_{L^{d+1}} \|f\|_{W^{1,2}}.\]
Therefore, $A:=FS_{1/(d+1)}$ is bounded relative to both $H$ and $H_0$.
Moreover, since $1/2=1/p_d'+1/(d+1)$,
\[ \|F g\|_{L^2} \le \|F\|_{L^{d+1}} \|g\|_{L^{p_d'}} \]
so that by definition of $X^*$,
\[ \|Af\|_{2} \le \|F\|_{L^{d+1}}\|f\|_{X^*}. \]
Hence, for any $I\subset \mathbb{R}\setminus( \mathcal{E}\cup\{0\})$,
\[ \sup_{0<\epsilon<1} \sup_{\lambda\in I} \|FS_{1/(d+1)} R_L(\lambda+i\epsilon) S_{1/(d+1)}
\overline{F}\|_{L^2\to L^2} \le C(I,L)\, \|F\|_{L^{d+1}}^2,\]
see~\eqref{bi3} and similarly with $R_0$. By Kato's
theorem~\cite{kato}, more precisely the local version of this
theorem as given by~\cite[Theorem~XIII.30]{RS4}, these properties
imply that $FS_{1/(d+1)}$ is smoothing relative to both $H$ and
$H_0$, and the constants involved only depend on $\|F\|_{d+1}$.
Using Kato's theory~\cite{kato},
\begin{equation*}
\int_{-\infty}^\infty||FS_{1/(d+1)}[e^{itH}E(I)f||_{L^2}^2\,dt\leq C(I,L)||f||_{L^2}^2||F||_{L^{d+1}}^2,
\end{equation*}
which is equivalent to the bound on the first term in \eqref{mai9}. For the second term we define
\[ A:=R^{-1/2}\mathbf{1}_{[|x|\le R]} S_1, \]
for any $R\geq 1$ and argue as before.

\section{Examples of admissible perturbations}\label{examples}

In this section we prove Proposition \ref{prop16}. We notice first that
\begin{equation*}
||Vf||_{B}\leq||V||_{Y}||f||_{B^\ast}
\end{equation*}
for any $V\in Y$ and $f\in B^\ast$. The constants $C$ in this section may depend on the exponent $q_0$ if $d=2$. Part (c) of Proposition \ref{prop16} is clear, directly from the definition of admissible perturbations.

For part (a) we prove the following lemma:

\newtheorem{lemma100}{Lemma}[section]
\begin{lemma100}\label{lemma100}
We have
\begin{equation}\label{tb102}
||\,|V|^{1/2}S_{-1/(d+1)}f||_{L^2}\leq
C||M_{q_0}(V)||_{L^{(d+1)/2}}^{1/2}||f||_{L^{p'_d}},\,f\in
L^{p'_d},
\end{equation}
\begin{equation}\label{tb300}
||\,|V|^{1/2}S_{-1}f||_{L^2}\leq
C||M_{q_0}(V)||_{Y}^{1/2}||f||_{B^\ast},\,f\in B^\ast,
\end{equation}
and
\begin{equation}\label{tb106}
||\,|V|^{1/2}S_{-1}f||_{L^2}\leq C||\,|V|\ast
K_{d,1/2}||_{Y}^{1/2}||f||_{B^\ast},\,f\in B^\ast.
\end{equation}
\end{lemma100}

\begin{proof}[Proof of Lemma \ref{lemma100}] We use the fact that for $\alpha\in\{1/(d+1),1\}$
\begin{equation*}
|S_{-\alpha}f(x)|\leq C|f|\ast W_\alpha(x),
\end{equation*}
where
\begin{equation*}
W_\alpha(x)=|y|^{-(d-\alpha)}\mathbf{1}_{\{|y|\leq
1\}}+|y|^{-(d+1)}\mathbf{1}_{\{|y|\geq 1\}}.
\end{equation*}
For any $s\in\mathbb{Z}^d$ let $Q_s$ denote the cube
$\{x:\sup_{i=1,\ldots,d}|x_i-s_i|\leq 1/2\}$. For \eqref{tb102},
using the Cauchy--Schwartz inequality and fractional integration
\begin{equation*}
\begin{split}
||\,|V|^{1/2}&S_{-1/(d+1)}f||_{L^2}^2\leq
C\sum_{s\in\mathbb{Z}^d}\int_{Q_s}|V(x)|[|f|\ast
W_{1/(d+1)}(x)]^2\,dx\\
&\leq C\sum_{s\in\mathbb{Z}^d}||V||_{L^{q_0}(Q_s)}\cdot||\,|f|\ast
W_{1/(d+1)}||_{L^{2q_0'}(Q_s)}^2\\
&\leq
C\sum_{s\in\mathbb{Z}^d}||V||_{L^{q_0}(Q_s)}\Big[\sum_{s'\in\mathbb{Z}^d}||\,(\mathbf{1}_{Q_{s'}}|f|)\ast
W_{1/(d+1)}||_{L^{2q_0'}(Q_s)}\Big]^2\\
&\leq
C\sum_{s\in\mathbb{Z}^d}||V||_{L^{q_0}(Q_s)}[\sum_{s'\in\mathbb{Z}^d}||f||_{L^{p'_d}(Q_{s'})} (1+|s-s'|)^{-d-1}]^2\\
&\leq
C\Big[\sum_{s\in\mathbb{Z}^d}||V||_{L^{q_0}(Q_s)}^{(d+1)/2}\Big]^{2/(d+1)}||f||_{L^{p'_d}}^2,
\end{split}
\end{equation*}
which gives \eqref{tb102}. The proof of \eqref{tb300} is similar:
\begin{equation*}
\begin{split}
||\,|V|^{1/2}&S_{-1}f||_{L^2}^2\leq
C\sum_{s\in\mathbb{Z}^d}\int_{Q_s}|V(x)|[|f|\ast
W_{1}(x)]^2\,dx\\
&\leq
C\sum_{s\in\mathbb{Z}^d}||V||_{L^{q_0}(Q_s)}\Big[\sum_{s'\in\mathbb{Z}^d}||\,(\mathbf{1}_{Q_{s'}}|f|)\ast
W_{1}||_{L^{2q_0'}(Q_s)}\Big]^2\\
&\leq
C\sum_{s\in\mathbb{Z}^d}||V||_{L^{q_0}(Q_s)}[\sum_{s'\in\mathbb{Z}^d}||f||_{L^2(Q_{s'})} (1+|s-s'|)^{-d-1}]^2\\
&\leq C\sum_{j=0}^\infty(2^j\sup_{s\in \mathbb{Z}^d\cap
D_j}||V||_{L^{q_0}(Q_s)})\cdot T_j
\end{split}
\end{equation*}
where, assuming $j$ fixed,
\begin{equation}\label{tb301}
\begin{split}
T_j&=2^{-j}\sum_{s\in\mathbb{Z}^d\cap
D_j}[\sum_{s'\in\mathbb{Z}^d}||f||_{L^2(Q_{s'})}
(1+|s-s'|)^{-d-1}]^2\\
&\leq C2^{-j}\sum_{s\in\mathbb{Z}^d\cap
D_j}\sum_{j'=0}^\infty2^{|j-j'|/10}[\sum_{s'\in\mathbb{Z}^d\cap D_{j'}}||f||_{L^2(Q_{s'})}
(1+|s-s'|)^{-d-1}]^2\\
&\leq C2^{-j}\sum_{j'=0}^\infty2^{|j-j'|/10}2^{-2|j-j'|}||f||_{L^2(D_{j'})}^2\leq C||f||_{B^\ast}^2.
\end{split}
\end{equation}
This completes the proof of \eqref{tb300}. To prove \eqref{tb106},
for any $s\in\mathbb{Z}^d$ let $\widetilde{Q}_s$ denote the cube
$\{x:\sup_{i=1,\ldots,d}|x_i-s_i|\leq 3/2\}$. We replace
fractional integration with the following local bound:
\begin{equation}\label{tb107}
||\,|V|^{1/2}[|f|\ast(|y|^{-(d-1)}\mathbf{1}_{\{|y|\leq 1\}})]||_{L^2(Q_s)}\leq C||\,|V|\ast K_{d,1/2}||_{L^\infty(\widetilde{Q}_s)}^{1/2}||f||_{L^2(\widetilde{Q}_s)}.
\end{equation}
This follows from \cite[Theorem 2.3]{KeSa}. Using the fact that
$||V||_{L^1(Q_s)}\leq C ||\,|V|\ast
K_{d,1/2}||_{L^\infty(\widetilde{Q}_s)}$, we have
\begin{equation*}
\begin{split}
||\,|V|^{1/2}&S_{-1}f||_{L^2}^2\leq
C\sum_{s\in\mathbb{Z}^d}\int_{Q_s}|V(x)|[|f|\ast
W_{1}(x)]^2\,dx\\
&\leq
C\sum_{s\in\mathbb{Z}^d}||\,|V|\ast K_{d,1/2}||_{L^\infty(\widetilde{Q}_s)}[\sum_{s'\in\mathbb{Z}^d}||f||_{L^2(Q_{s'})} (1+|s-s'|)^{-d-1}]^2\\
&\leq C\sum_{j=0}^\infty(2^j\sup_{s\in \mathbb{Z}^d\cap
D_j}||\,|V|\ast K_{d,1/2}||_{L^\infty(\widetilde{Q}_s)})\cdot T_j,
\end{split}
\end{equation*}
where $T_j$ is as above. The bound \eqref{tb301} completes the proof of the lemma.
\end{proof}

We return to the proof of Proposition \ref{prop16}. Let
$\mathcal{N}_1(V)=||M_{q_0}(V)||_{L^{(d+1)/2}}$, $
\mathcal{N}_2(V)=||M_{q_0}(V)||_{Y}$, and $
\mathcal{N}_3(V)=||\,|V|\ast K_{d,1/2}||_{Y}$. It follows from
Lemma \ref{lemma100} that
\begin{equation}\label{tb310}
||\,|V|^{1/2}u||_{L^2}\leq
C\min_{i\in\{1,2,3\}}\mathcal{N}_i(V)^{1/2}||u||_{X^\ast}
\end{equation}
for any $u\in X^\ast$. For potentials $V$ as in \eqref{tb99}, \eqref{tb991}, or \eqref{tb100} and $u\in X^\ast$ we define the distribution $L_Vu$ by the formula
\begin{equation}\label{tb311}
\langle L_Vu,\phi\rangle:=\langle |V|^{1/2}u,|V|^{1/2}\mathrm{sign}(V)\phi\rangle=\int_{\mathbb{R}^d}Vu\overline{\phi}\,dx.
\end{equation}
The distribution $L_Vu$ is well defined, in view of \eqref{tb310}. Using \eqref{tb511} and \eqref{tb310}
\begin{equation}\label{tb421}
||L_Vu||_X\leq
C\min_{i\in\{1,2,3\}}\mathcal{N}_i(V)||u||_{X^\ast},\,u\in X^\ast,
\end{equation}
thus $L_V\in\mathcal{L}(X^\ast,X)$. The identity \eqref{tb311}
also shows that $L$ is symmetric in the sense on \eqref{tb61}.

Next, we verify \eqref{tb6}. Let $\varphi:\mathbb{R}^d\to[0,C]$
denote a smooth function supported in the ball $\{x:|x|\leq 1\}$
with $\int_{\mathbb{R}^d}\varphi\,dx=1$, and
$\chi:\mathbb{R}^d\to[0,1]$ a smooth function supported in the
ball $\{x:|x|\leq 2\}$ and equal to $1$ in the ball $\{x:|x|\leq
1\}$. For $\varepsilon\in(0,1]$ and $r\in[1,\infty)$, let
$\varphi_\varepsilon(x)=\varepsilon^{-d}\varphi(x/\varepsilon)$
and $\chi_r(x)=\chi(x/r)$. For integers $n\geq 1$ let
\begin{equation*}
V_{n}(x)=\chi_n(x)(V\ast\varphi_{1/n})(x).
\end{equation*}
We will show that
\begin{equation}\label{tb420}
\begin{cases}
\text{ if }\mathcal{N}_1(V)<\infty&\text{ then
}\lim_{n\to\infty}\mathcal{N}_1(V-V_n)=0;\\
\text{ if }\mathcal{N}_2(V)<\infty&\text{ then
}\lim_{n\to\infty}\mathcal{N}_2(V-V_n)=0;\\
\text{ if }\lim_{\delta\to 0}||\,|V|\ast
K_{d,\delta}||_{Y}=0&\text{ then
}\lim_{n\to\infty}\mathcal{N}_3(V-V_n)=0.
\end{cases}
\end{equation}
Assuming \eqref{tb420}, the proof of \eqref{tb6} is easy. For
$i\in\{1,2,3\}$, given $\varepsilon$ as in \eqref{tb6}, we fix
$n=n(\varepsilon)$ with the property that
$\mathcal{N}_i(V-V_n)\leq(\varepsilon/C)^{1/2}$, where $C$ is the
constant in \eqref{tb421}. Using \eqref{tb421}
\begin{equation}\label{tb112}
\begin{split}
||\mu_{N,\gamma}L_Vu||_X&\leq ||L_{V-V_n}(\mu_{N,\gamma}u)||_X+||\mu_{N,\gamma}V_nu||_X\\
&\leq\varepsilon||\mu_{N,\gamma}u||_{X^\ast}+C||\mu_{N,\gamma}V_nu||_B\\
&\leq\varepsilon||\mu_{N,\gamma}u||_{X^\ast}+C_{V,N,\varepsilon}||u\mathbf{1}_{\{|x|\leq 2n\}}||_{L^2},
\end{split}
\end{equation}
as desired. It remains to verify \eqref{tb420}. The first two
limits in \eqref{tb420} are straightforward. For the last limit
fix $\varepsilon>0$. By the definition of the space $Y$, there is
$n_{\varepsilon,V}$ with the property that
\begin{equation*}
||[|V|\ast K_{d,1/2}]\mathbf{1}_{\{|x|\geq
n_{\varepsilon,V}\}}||_Y\leq\varepsilon/C.
\end{equation*}
Then
\begin{equation}\label{tb116}
||[|V-V_n|\ast K_{d,1/2}]\mathbf{1}_{\{|x|\geq
n_{\varepsilon,V}\}}||_Y\leq\varepsilon/3
\end{equation}
for any integer $n\geq 1$. The condition \eqref{tb100} shows that
there is $\delta_{\varepsilon,V}$ with the property that
\begin{equation}\label{117}
||[|V-V_n|\ast K_{d,\delta_{\varepsilon,V}}]\mathbf{1}_{\{|x|\leq
n_{\varepsilon,V}\}}||_Y\leq C||V\ast
K_{d,\delta_{\varepsilon,V}}||_Y\leq\varepsilon/3
\end{equation}
for any integer $n\geq 1$. Finally, notice that the kernel
$K_{d,1/2}-K_{d,\delta_{\varepsilon,V}}$ is bounded. Since $V\in
L^1_{\mathrm{loc}}(\mathbb{R}^d)$,
$\lim_{n\to\infty}[V_n-V]\mathbf{1}_{\{|x|\leq
n_{\varepsilon,V}+1\}}=0$ in $L^1$, so \eqref{tb420} follows.

To verify condition~(3) let $J=1$ and $A_1u:=|V|^{1/2}u$,
$B_1u:=|V|^{1/2}\mathrm{sign}(V)u$ with domains
\[ \mathrm{Domain}(A_1)=\mathrm{Domain}(B_1):=\{ f\in L^2: |V|^{1/2} f \in L^2 \}. \]
It follows from \eqref{tb310} that $A_1,B_1\in\mathcal{L}(X^\ast,L^2)$ and $W^{1,2}\subseteq\mathrm{Domain}(A_1)=\mathrm{Domain}(B_1)$. Also,
it follows from Fatou's lemma that $A_1, B_1$ are closed on this
domain. It remains to verify the identity \eqref{eq:Lsplit}. The
identity is clear  for $\phi,\psi\in\mathcal{S}(\mathbb{R}^d)$, in
view of \eqref{tb311}. For $\phi,\psi\in X^\ast$, we define the sequences $\phi_n$ and $\psi_n$ as in the proof of Lemma \ref{lemma9.1}. In view of \eqref{tb660}, which was proved using only conditions (1) and (2) in the definition of admissible perturbations, it suffices to prove that
\begin{equation}\label{tb670}
\lim_{n\to\infty}\langle A_1\phi_n,B_1\psi_n\rangle=\langle A_1\phi,B_1\psi\rangle.
\end{equation}
We may assume that $||\phi||_{X^\ast}=||\psi||_{X^\ast}=1$. Given $\varepsilon_0>0$, we fix $n_0$ with the property that $\min_{i\in\{1,2,3\}}\mathcal{N}_i(V-V_{n_0})\leq\varepsilon$ (using \eqref{tb420}). Using \eqref{tb310}
\begin{equation*}
\begin{split}
&|\langle A_1\phi,B_1\psi\rangle-\int_{\mathbb{R}^d}V_{n_0}\phi\overline{\psi}\,dx|\leq C\varepsilon;\\
&|\langle A_1\phi_n,B_1\psi_n\rangle-\int_{\mathbb{R}^d}V_{n_0}\phi_n\overline{\psi}_n\,dx|\leq C\varepsilon.
\end{split}
\end{equation*}
The limit \eqref{tb670} follows since $\lim_{n\to\infty}\phi_n\mathbf{1}_{\{|x|\leq 2n_0\}}=\phi\mathbf{1}_{\{|x|\leq 2n_0\}}$ in $L^2$ and $\lim_{n\to\infty}\psi_n\mathbf{1}_{\{|x|\leq 2n_0\}}=\psi\mathbf{1}_{\{|x|\leq 2n_0\}}$ in $L^2$.
\medskip

We now prove part (b) of the proposition. Using part (c), we may
assume that $\vec{a}=(0,\ldots,0,a)$, so
\begin{equation*}
\vec{a}\cdot\nabla-\nabla\cdot\overline{\vec{a}}=a\partial_{x_d}-\partial_{x_d}\overline{a}.
\end{equation*}
We are looking to define the distribution $L_{a}$ by the formula
\begin{equation}\label{tb605}
\langle
L_{a}u,\phi\rangle:=\langle\omega\partial_{x_d}u,\omega^{-1}\overline{a}\phi\rangle+\langle\omega^{-1}\overline{a}u,\omega\partial_{x_d}\phi\rangle,
\end{equation}
for any $u\in X^\ast$ and $\phi\in\mathcal{S}(\mathbb{R}^d)$. Here
\begin{equation}\label{tb606}
\omega=\sum_{j=0}^\infty2^{-j/2}\omega_j\mathbf{1}_{D_j},
\end{equation}
where $\omega_j>0$ are real numbers that will be fixed depending
on the function $a$. The distribution $L_a$ in \eqref{tb605} is
well defined if
\begin{equation*}
||\omega\partial_{x_d}||_{X^\ast\to L^2}<\infty
\end{equation*}
and
\begin{equation*}
||\omega^{-1}a||_{X^\ast\to L^2}<\infty.
\end{equation*}
Since $X^\ast\subseteq S_{-1}B^\ast$, we have
\begin{equation}\label{tb609}
||\omega\partial_{x_d}||_{X^\ast\to L^2}\leq C\Big[\sum_{j=0}^\infty\omega_j^2\Big]^{1/2}.
\end{equation}
Assume that the sequence $\omega_j$ is chosen in such a way that
\begin{equation}\label{tb607}
C^{-1}\omega_j\leq \omega_{j+1}\leq C\omega_j
\end{equation}
for any integer $j\geq 0$. Then, using \eqref{tb102} we have
\begin{equation}\label{tb610}
\begin{split}
||\omega^{-1}&au||_{L^2}\leq C||M_{q_0}(\omega^{-2}a^2)||_{L^{(d+1)/2}}^{1/2}||u||_{X^\ast}\\
&\leq C||u||_{X^\ast}\Big[\sum_{j=0}^\infty\int_{D_j}[M_{q_0}(\omega^{-2}a^2)]^{(d+1)/2}\,dx\Big]^{1/(d+1)}\\
&\leq
C||u||_{X^\ast}\Big[\sum_{j=0}^\infty\omega_j^{-(d+1)}[2^{j/2}||M_{2q_0}(a)||_{L^{d+1}(D_j)}]^{d+1}\Big]^{1/(d+1)}.
\end{split}
\end{equation}
Using \eqref{tb300} we have
\begin{equation}\label{tb611}
\begin{split}
||\omega^{-1}&au||_{L^2}\leq C||M_{q_0}(\omega^{-2}a^2)||_{Y}^{1/2}||u||_{X^\ast}\\
&\leq C||u||_{X^\ast}\Big[\sum_{j=0}^\infty2^j||M_{q_0}(\omega^{-2}a^2)||_{L^\infty}(D_j)\Big]^{1/2}\\
&\leq
C||u||_{X^\ast}\Big[\sum_{j=0}^\infty\omega_j^{-2}[2^j||M_{2q_0}(a)||_{L^\infty(D_j)}]^{2}\Big]^{1/2}.
\end{split}
\end{equation}
Using \eqref{tb106} we have
\begin{equation}\label{tb612}
\begin{split}
||\omega^{-1}&au||_{L^2}\leq C||\omega^{-2}|a|^2\ast K_{d,1/2}||_{Y}^{1/2}||u||_{X^\ast}\\
&\leq
C||u||_{X^\ast}\Big[\sum_{j=0}^\infty\omega_j^{-2}[2^j||(|a|^2\ast
K_{d,1/2})^{1/2}||_{L^\infty(D_j)}]^2\Big]^{1/2}.
\end{split}
\end{equation}

To deal with potentials $a$ as in \eqref{tb98}, we would like to
fix
\begin{equation*}
\omega_j=C_a[2^{j/2}||M_{2q_0}(a)||_{L^{d+1}(D_j)}]^{(d+1)/(d+3)},
\end{equation*}
in order to optimize \eqref{tb609} and \eqref{tb610}. This is not
possible because of the restriction \eqref{tb607}. To avoid this
problem, let
\begin{equation*}
\theta_j=\sum_{j'=0}^\infty
2^{j'/2}||M_{2q_0}(a)||_{L^{d+1}(D_{j'})}2^{-|j-j'|}\text{ and
}\omega_j=[\theta_j]^{\frac{d+1}{d+3}}/
\Big(\sum_{j=0}^\infty[\theta_j]^{p_d}\Big)^{\frac{d-1}{4(d+1)}}.
\end{equation*}
Clearly, $2^{j/2}||M_{2q_0}(a)||_{L^{d+1}(D_{j})}\leq \theta_j$ and
\eqref{tb607} holds. By \eqref{tb609} and \eqref{tb610}
\begin{equation}\label{tb620}
||\omega\partial_{x_d}u||_{L^2}+||\omega^{-1}au||_{L^2}\leq
C||u||_{X^\ast}\Big[\sum_{j=0}^\infty[2^{j/2}||M_{2q_0}(a)||_{L^{d+1}(D_j)}]^{p_d}\Big]^{1/(2p_d)}.
\end{equation}

Similarly, to deal with potentials as in \eqref{tb981}, we let
\begin{equation*}
\theta_j=\sum_{j'=0}^\infty
2^{j'}||M_{2q_0}(a)||_{L^\infty(D_{j'})}2^{-|j-j'|}\text{ and
}\omega_j=[\theta_j]^{1/2},
\end{equation*}
and it follows from \eqref{tb609} and \eqref{tb611} that
\begin{equation}\label{tb621}
||\omega\partial_{x_d}u||_{L^2}+||\omega^{-1}au||_{L^2}\leq
C||u||_{X^\ast}\Big[\sum_{j=0}^\infty2^{j}||M_{2q_0}(a)||_{L^\infty(D_{j})}\Big]^{1/2}.
\end{equation}
Finally, to deal with potentials as in \eqref{tb100}, we let
\begin{equation*}
\theta_j=\sum_{j'=0}^\infty 2^{j'}||(|a|^2\ast
K_{d,1/2})^{1/2}||_{L^\infty(D_{j'})}2^{-|j-j'|}\text{ and
}\omega_j=[\theta_j]^{1/2},
\end{equation*}
and it follows from \eqref{tb609} and \eqref{tb612} that
\begin{equation}\label{tb622}
||\omega\partial_{x_d}u||_{L^2}+||\omega^{-1}au||_{L^2}\leq
C||u||_{X^\ast}\Big[\sum_{j=0}^\infty2^{j}||(|a|^2\ast
K_{d,1/2})^{1/2}||_{L^\infty(D_{j})}\Big]^{1/2}.
\end{equation}

It follows from \eqref{tb620}, \eqref{tb621} and \eqref{tb622}
that the distribution $L_a$ in \eqref{tb605} is well defined. In
fact,
\begin{equation}\label{tb630}
\langle L_au,\phi\rangle=\int_{\mathbb{R}^d}a(\partial_{x_d}u)
\overline\phi+\overline{a}u(\partial_{x_d}\overline\phi)\,dx,
\end{equation}
for any $u\in X^\ast$ and $\phi\in\mathcal{S}(\mathbb{R}^d)$,
where the integral converges absolutely. Let $
\mathcal{N}'_1(a)=\Big[\sum_{j=0}^\infty(2^{j/2}||M_{2q_0}(a)||_{L^{d+1}(D_j)})^{p_d}\Big]^{1/p_d}$,
$\mathcal{N}'_2(a)=||M_{2q_0}(a)||_{Y}$, and
$\mathcal{N}'_3=||(|a|^2\ast K_{d,1/2})^{1/2}||_{Y}$. Using
\eqref{tb511}, it follows from \eqref{tb605}, \eqref{tb620},
\eqref{tb621}, and \eqref{tb622} that
\begin{equation}\label{tb652}
||L_au||_{X}\leq
C||u||_{X^\ast}\min_{i\in\{1,2,3\}}\mathcal{N}'_i(a),
\end{equation}
Thus $L_a\in\mathcal{L}(X^\ast,X)$. It follows easily from
\eqref{tb605} that $L_a$ is symmetric, in the sense of
\eqref{tb61}.

To prove condition $(2)$ in the definition of admissible
perturbations, let
\begin{equation*}
a_n(x)=\chi_n(x)(a\ast\varphi_{1/n})(x),
\end{equation*}
where $\chi$ and $\varphi$ are as before. As in the proof of
\eqref{tb420}, it follows that
\begin{equation*}
\begin{cases}
\text{ if }\mathcal{N}'_1(a)<\infty&\text{ then
}\lim_{n\to\infty}\mathcal{N}'_1(a-a_n)=0;\\
\text{ if }\mathcal{N}'_2(a)<\infty&\text{ then
}\lim_{n\to\infty}\mathcal{N}'_2(a-a_n)=0;\\
\text{ if }\lim_{\delta\to 0}||\,[|a|^2\ast
K_{d,\delta}]^{1/2}||_{Y}=0&\text{ then
}\lim_{n\to\infty}\mathcal{N}'_3(a-a_n)=0.
\end{cases}
\end{equation*}
The identity \eqref{tb630} shows that $L_a=L_{a_n}+L_{a-a_n}$. An
estimate similar to \eqref{tb112} shows that it suffices to prove
that for any $N\geq 0$
\begin{equation}\label{tb640}
||\mu_{N,\gamma}L_{a}u||_X\leq
C_N\min_{i\in\{1,2,3\}}\mathcal{N}'_i(a)||\mu_{N,\gamma}u||_{X^\ast},
\end{equation}
for any $u\in X^\ast$ and any $\gamma\in(0,1]$. With the notation in Lemma \ref{lemma6},
it is easy to check that
\begin{equation*}
\mu_{N,\gamma}L_{a}u=L_a(\mu_{N,\gamma}u)-b_d(a-\overline{a})(\mu_{N,\gamma}u).
\end{equation*}
The bound for the first term follows directly from \eqref{tb652}.
For the second term, we define $\omega=\omega(a)=\omega(\overline{a})$
as before, and use \eqref{tb511}, \eqref{tb620}, \eqref{tb621}, and \eqref{tb622}:
\begin{equation*}
\begin{split}
||b_d(a-\overline{a})&(\mu_{N,\gamma}u)||_{X}\leq C\sup_{\phi\in\mathcal{S}(\mathbb{R}^d),\,||\phi||_{X^\ast}=1}|\langle \omega^{-1}(a-\overline{a})(\mu_{N,\gamma}u),\omega b_d\phi\rangle|\\
&\leq C||\mu_{N,\gamma}u||_{X^\ast}\min_{i\in\{1,2,3\}}\mathcal{N}'_i(a)^{1/2}\sup_{\phi\in\mathcal{S}(\mathbb{R}^d),\,||\phi||_{X^\ast}=1}||\omega b_d\phi||_{L^2}\\
&\leq C||\mu_{N,\gamma}u||_{X^\ast}\min_{i\in\{1,2,3\}}\mathcal{N}'_i(a),
\end{split}
\end{equation*}
using \eqref{tb651}. This completes the proof of \eqref{tb6}.

For condition (3) in the definition of admissible perturbations, let $J=2$ and
\begin{equation*}
A_1u=B_2u:=\omega^{-1}\overline{a}u,\text{ and }B_1u=A_2u:=\omega\partial_{x_d}u,
\end{equation*}
where $\omega$ is defined as before. As domains we again choose
the natural ones, i.e.,
\begin{align*}
\mathrm{Domain}(A_1) &=\mathrm{Domain}(B_2):=\{ f\in L^2:
\omega^{-1}\overline{a}f \in L^2 \} \\
\mathrm{Domain}(B_1) &=\mathrm{Domain}(A_2):=\{ f\in L^2:
\omega\partial_{x_d}f \in L^2 \}.
\end{align*}
It is again easy to see that these domains make $A_1, A_2,B_1,B_2$
closed.\footnote{Strictly speaking, one should smooth out
$\mathbf{1}_{D_j}$ in \eqref{tb606} here, but we ignore this
issue.}
 It follows from
\eqref{tb620}, \eqref{tb621}, and \eqref{tb622} that
$A_1,\,B_1,\,A_2,\,B_2\in\mathcal{L}(X^\ast,L^2)$. To verify the
identity \eqref{eq:Lsplit}, we notice first that the identity
holds if $\phi,\psi\in\mathcal{S}(\mathbb{R}^d)$, in view of the
definition \eqref{tb605}. The proof for $\phi,\psi\in X^\ast$ then
follows by the same limiting argument as in part (a), see the
proof of \eqref{tb670}.

\end{document}